%% file: BD-mTSP.tex
\documentclass[times,5p,twocolumn,authoryear]{elsarticle}

\journal{computers \& operations research}
\input{wolfpack-mTSP}

\begin{document}
\begin{frontmatter}
\title{Balanced dynamic multiple travelling salesmen:\\ algorithms and continuous approximations}

\author{Wolfgang Garn\fnref{w.garn@surrey.ac.uk}}
\address{Department of Business Transformation, University of Surrey, United Kingdom}

\input{intro-abstract}
\begin{keyword}
	travelling salesman \sep dynamic routing \sep mTSP \sep
    transportation \sep continuous approximation model
\end{keyword}
\end{frontmatter}

\input{sec-Introduction}
\input{sec-Dynamics}
\input{sec-Heuristics}

\input{sec-real-world-apps}
\input{sec-CAM}

\input{sec-Conclusion}

\section{Acknowledgement} 
My thanks goes to the inspiring and creative Hyster-Yale UK team, 
whose work on the next generation of intelligent autonomous forklift-trucks motivated me to look for mTSP solutions. 
Mark's vision of forklift-trucks bargaining for the next job inspired the closest vehicle mTSP heuristic, 
which is suitable for the online-mTSP, i.e. a mTSP without knowing the pallet tasks (or customer locations) in advance. 
Chi (Ezeh) encouraged me to keep a wholistic view that reached from strategic mTSP formulations to vehicles' Artificial Intelligence. 
Chris drew my attention to concurrency issues, that vehicles encounter when using the same path. 
There were many more inspirations, which demonstrated the importance of the online-mTSP and its solution approaches in practical applications. 

I am grateful for the anonymous reviewer of my work \citep{garn2020closed};
who provided many interesting pointers to valuable research articles,
which was the motivation for this work.
I would like to thank the reviewers of this article for encouraging me to consider real-world applications,
which added another dimension to this work.
Overall, their feedback greatly improved the quality of this work.
\section{Appendix}
The appendix provides five tables detailing the obtained distances for four sets of test instances: set 1, set 2, set X-U and set XXL.
These were analysed in respect to $m$-absolute dynamics (i.e. based on vehicles)  and relative dynamics (i.e. based on customers). 
The considered $m$-absolute dynamics $D_a^m$ use $\set{\frac{1}{2}, 1, \frac{3}{2}, 2, 4, 8}$,
and the relative dynamics are $D_r \in \set{2\%, 5\%, 7\%, 10\%, 20\%, 30\%, 100\%}$.
Corresponding absolute dynamics $D_a$ and 
the number of vehicles $m$ are shown in the tables.
Table \ref{tab:set1-abs} shows the $m$-absolute dynamics for set 1, 
which includes several new micro mTSP instances and some derived from ``classic'' TSP instances.
The corresponding relative dynamic results are given within the text in Table \ref{tab:set1-rel}.
Table \ref{tab:set2-rel} and \ref{tab:set2-abs} show $D_a^m$ and $D_r$ for set 2,
which were derived from the E\&K TSP instances.
Table \ref{tab:setX-rel} and \ref{tab:setX-abs} show $D_a^m$ and $D_r$ is based on the CVRP X test instance (specifically those with unitary demand),
which were obtained from \cite{uchoa2017new}.
Note that $D_r=100\%$ means that BD-CVRP is identical to the B-CVRP. 
The B-CVRP can be compared to the CVRP.
Table \ref{tab:setXXL-abs} contains extra-large test instances with more than 1,000 customers.
The modified test instances origin from \cite{arnold2019efficiently}.

\input{tab-set1-abs}

\input{tab-set2-rel}
\input{tab-set2-abs}

\input{tab-setx-rel}
\input{tab-setx-abs}
\input{tab-setXXL-abs}
\bibliographystyle{elsarticle-harv} 
\bibliography{mTSP}
\end{document}

%% file: wolfpack-mTSP.tex
\usepackage{lineno} 
\modulolinenumbers[1] 

\usepackage{amsmath,amsfonts,amssymb,amsthm,mathrsfs,graphicx,booktabs} 
\usepackage[export]{adjustbox} 
\usepackage{bigstrut,multirow,rotating}
\usepackage{hyperref,colortbl,algorithm,algorithm,algpseudocode}
\usepackage{lmodern}

\theoremstyle{definition}
\newtheorem{definition}{Definition}[section]

\def\real{\mathbb{R}} 
 \def\setN{\mathbb{N}} \def\setB{\mathbb{B}}

\DeclareMathOperator*{\argmin}{arg\,min}
\newcommand{\set}[1]{\left\{#1\right\}}
\newcommand{\setn}{\{1,2,\dots,n\}}
\newcommand{\sett}[1]{\{1,2,\dots,#1\}}
\newcommand{\seq}[1]{\left\langle#1\right\rangle}
\newcommand{\floor}[1]{\left\lfloor #1 \right\rfloor}%
\newcommand{\ceil}[1]{\left\lceil #1 \right\rceil}%
\newcommand{\round}[1]{\left\lfloor #1 \right\rceil}
\newcommand{\norm}[1]{\left\lVert#1\right\rVert}
\newcommand{\mat}[1]{\begin{bmatrix}#1\end{bmatrix}}
\newcommand{\equ}[1]{\begin{equation}#1\end{equation}}
\newcommand{\lequ}[2]{\begin{equation}#2\label{eq:#1}\end{equation}} 
\newtheorem{example}{Example}[section]
\newtheorem{theorem}{Theorem}
\newtheorem{proposition}{Proposition}

\newcommand{\obj}[2][max]{& \text{#1} & & #2 \\}  
  
\newcommand{\sto}[1]{& & & #1 \\}  
\newcommand{\LPs} [2]{\begin{aligned} #1 & \text{subject to} & & \\ #2 \end{aligned}}
\newcommand{\LP}  [2]{\begin{equation}\LPs{#1}{#2} \end{equation}}
\newcommand{\mcomment}[1]{\text{\textcolor{cyan}{~~$\triangleright$\textit{~#1}}}} 
\newcommand{\ie}{\Rightarrow}
\newcommand{\tp}[1]{\ensuremath{#1^\intercal}} 

\newcommand{\code}[1]{\texttt{#1}}

\DeclareMathOperator{\arctantwo}{arctan2}
\newcommand {\ms}[1]{\ensuremath{\overline{#1}\pm \sigma_{#1}}}

%% file: intro-abstract.tex
\begin{abstract}
Dynamic routing occurs when customers are not known in advance, e.g. for real-time routing.
Two heuristics are proposed that solve the balanced dynamic multiple travelling salesmen problem (BD-mTSP). 
These heuristics represent operational (tactical) tools for dynamic (online, real-time) routing.
Several types and scopes of dynamics are proposed.
Particular attention is given to sequential dynamics.
The balanced dynamic closest vehicle heuristic (BD-CVH) and
the balanced dynamic assignment vehicle heuristic (BD-AVH) are applied to this type of dynamics.
The algorithms are applied to a wide range of test instances.
Taxi services and palette transfers in warehouses demonstrate how to use the BD-mTSP algorithms in real-world scenarios.

Continuous approximation models for the BD-mTSP's are derived
and serve as strategic tools for dynamic routing.
The models express route lengths using vehicles, customers, and dynamic scopes without the need of running an algorithm.
A machine learning approach was used to obtain regression models.
The mean absolute percentage error of two of these models is below 3\%.
\end{abstract}

%% file: sec-Introduction.tex
\section{Introduction}\label{sec:Introduction} 
In several applications such as
taxi-services \citep{caramia2002routing,fabri2006dynamic}, 
emergency-dispatches \citep{yang2005online} and 
warehouse-picking \citep{smolic2009time}
the next location is unknown.
In such cases the multiple Travelling Salesman Problem (mTSP) becomes dynamic.
That means during the ``building'' of the route new nodes are revealed.
Hence, there is either an implicit or explicit time dimension involved.

In the \textit{multiple Travelling Salesman Problem} (mTSP) 
$m$ salesmen (vehicles) visit $n$ customers (nodes).
Each customer will be visited once by one and only one salesman.
The salesmen start and end their routes at a single depot (source node).
The objective is to minimise the total distance travelled by the salesmen.
The mTSP is called \textit{balanced} 
when each vehicle visits ``approximately'' the same number of customers.
This can be achieved by introducing an upper and lower limit of customers a salesman must visit.
\cite{gouveia2010vehicle} discuss these capacity bounds in the context of the Vehicle Routing Problem.
\cite{martinez2015transformations} introduce a general approach that transforms node-balanced routing problems 
into generalised balanced TSPs. This allows the reduction of arcs in the underlying graph.
\cite{bektacs2019balanced} solve balanced vehicle routing using a polyhedral analysis and a branch-and-cut algorithm.
An interesting integer program formulation of the balanced-mTSP is given in \cite{kara2006integer}.
I adapted this program in \cite{garn2020closed} for the balanced-mTSP 
and compared it to two heuristics, one of them is the pre-cursor for this work.
That heuristic focused on the static balanced mTSP, 
whereas the heuristics in this paper deal with dynamics.
The earlier work offers a continuous approximation formula for the static balanced-mTSP 
using a simple regression approach,
whilst the work here introduces a framework, which has not been used for continuous approximation models previously. 

The aim of this work is to provide a continuous approximation model (CAM) for the balanced-dynamic mTSP (BD-mTSP)
proposing a Machine Learning approach.

Hence, the rest of the paper is organised as follows.
Section \ref{sec:Dynamics} classifies types of dynamics.
Sequential absolute dynamics will be used to develop
two BD-mTSP heuristics (BD-CVH and BD-AVH in Section \ref{sec:Heuristics}).
They are compared against each other and an exact balanced-static (BS) mTSP method (Section \ref{ssec:heuristics-test-instances}). 
Here, a number of classic TSP/mTSP test instances are used.
It must be stressed that the comparison is intended 
to derive the difference between dynamic and static mTSPs rather than whether one is better than the other one.
This will give first insights of how much balanced-static and balanced-dynamic mTSPs differ from each other.
Section \ref{sec:real-world-apps} explains 
how to apply the BD-mTSP to real-world instances such as taxi services and transfer of pallets in warehouses.
Continuous approximation models for the static mTSP are reviewed in Section \ref{sec:CAM}.
Most approaches focus on uniformly distributed nodes in the Euclidean plane,
which will be adapted for the dynamic case.
To the best of the author's knowledge this is the first CAM model for the balanced-dynamic case.
Moreover, this appears to be that for the first time that a structured Machine Learning approach
was used to derive a CAM relationship.
The here introduced BD-mTSP continuous approximation models describe the total distance depending on 
the number of vehicles, customers, and sequential-dynamics. 

%% file: sec-Dynamics.tex
\section{Dynamics \label{sec:Dynamics}}
\subsection{Evolution and Quality}

\cite{pillac2013review} reviewed dynamic vehicle routing problems, 
which is also relevant to the mTSP.
They suggested to differentiate between information evolution and quality.
I will propose definitions to make their ideas more tangible.
\textit{Information evolution} refers to changes in the data available to the planner such as new customer requests.
Let $C=\seq{c_1,\ldots,c_n}=\seq{1,\ldots,n}$ be the sequence of customers changing over an evolving time horizon $t$.
$n$ denotes the finite maximum number of customers for all $t$.
Assume that $t \in T=\seq{t_1,\ldots,t_n}$ is a discrete time event where an information evolution occurred.
Without restricting generality let $C$ be ordered over time, 
i.e. if $i<j$ then $c_i$ information evolution happened before or at the same time as $c_j$.
\begin{definition}{\textit{Node information evolution}}
	is the process of associating time $t_i \in T$ to nodes $c_i \in C$. 
	$T$ is a monotonic increasing sequence.
	$C$ is ordered such that each $t_i$ is mapped to $c_i$.
\end{definition}
Hence, the information evolution is the ``visibility'' of customers at a certain point in time.
The above definition automatically ensures that related information 
such as the location or distances between customers becomes time dependent.
That means, the customer locations $X$ can be deterministic but revealed over time.
In the Euclidean plane $X=\mat{x&y}$ with coordinates $x=\mat{x_1,\ldots,x_n}$ and $y=\mat{y_1,\ldots,y_n}$.
The above definition means that customers and locations are ordered according to a time dimension.
A special case, which will be used later, is that all time steps are discrete and equidistant.

The sequence of vehicles (fleet) $V=\seq{v_1,\ldots,v_n}=\seq{1,\ldots,m}$ and related information can evolve over time as well.
Hence, \textit{vehicle information evolution} is defined similarly.
Later sections explaining heuristics and CAMs will assume that the fleets existence is permanent.

Information quality is related to uncertainty in the input data.

\begin{definition}{\textit{Node information quality}}
	is the probability distribution associated to the positions of nodes
	or distances between nodes.
\end{definition}
An example for requiring an exact (deterministic) location is the picking of an item.
The eventual location of an object in a drop-off zone is an example for an uncertain (stochastic) position.

Table \ref{tab:evolution-quality} shows the two dimensions and its related classes.
\begin{table}[htbp]	\centering	\caption{Information evolution and quality.}
	\begin{tabular}{cc|c|c}		\toprule
		&      & \multicolumn{2}{c}{\textbf{quality}} \\
		& \textbf{input} & deterministic & stochastic \\
		\midrule \rule{0pt}{17pt}
		\multirow{2}[0]{*}[3mm]{\begin{sideways}\textbf{evolution}\end{sideways}} & static & SD   & SS \\
		\cmidrule{2-4}  \rule{0pt}{17pt}         & dynamic & DD   & DS \\		\bottomrule
	\end{tabular} \label{tab:evolution-quality}
\end{table}
The ``classic'' mTSP falls into the SD class, i.e. static (no) information evolution and deterministic information quality.
The exact balanced mTSP introduced in \cite{garn2020closed} is in this class.
Even mTSPs with time-windows are in this class. 
In the context of routing - stochastic means that information such as customer locations follows a known random distribution.
Other typical stochastic VRP factors are demand, times, and pick-ups.
On a strategic level all the static continuous approximation models (see Section \ref{sec:CAM}) fall into the SS class.
Any VRP or mTSP implementation that is supplied with data based on a random distribution is called stochastic.
The balanced dynamic closest vehicle heuristic (BD-CVH) and balanced dynamic assignment vehicle heuristic (BD-AVH) are able to operate in the dynamic-deterministic (DD) class 
where information evolution is dynamic 
and information quality is deterministic (all locations are known).
The DD class means for the mTSP that some of the $n$ customers will not be known in advance,
but at some time during the execution of the route.
The dynamic-stochastic (DS) class implies for the mTSP that new customers are revealed during the building of the route,
and their locations or distances are randomly distributed.
Note that some of the customer locations could be known in advance.
The knowledge of the underlying random distribution allows algorithms to anticipate (predict) the occurrence of new locations.
The continuous approximation models given later are prime examples for the DS-class.
The BD-CVH and BD-AVH operate in this class as well.

\cite{pillac2013review} mention dynamic VRP solution methods and divide
them into continuous and periodic approaches.
Here, continuous, implies as soon as new information becomes available re-optimisation is performed.
An example for continuous re-optimisation is provided by \cite{gendreau1999parallel}.
They use parallel tabu search to find solutions for the dynamic VRP.
Their premises are that new information is an event that triggers a re-optimisation.
This allowed them to adapt a static mTSP implementation in a dynamic context.
However, the heuristics focus on accommodating minor changes.
For comparison purposes they proposed a few algorithms. 
Their ``insertion'' algorithm adds a new customer to the planned routes such that the added cost is minimised. 
A ``rebuild'' algorithm and a few more adaptive tabu search methods are briefly mentioned as options.
Hence, if most information is known in advance and only minor adaptation are expected these algorithms are a good option.
This type of dynamics will be called \textit{random node insertion}.

Genetic algorithms (GA) are a popular approach, and so they find their presence in the continuous D-mTSP world. 
\cite{cheung2008dynamic} used a GA to update any changes to existing routes.
Again, this method is based on the assumption of having most of the information available.
Similarly, \cite{haghani2005dynamic} developed a GA to solve a D-VRP with time-dependent travel times. 
D-mTSP based on periodic re-optimisation fall back on static mTSP solution methods.
Whenever, new information is available all data is assumed to be static and deterministic 
and a classic mTSP solution procedure is executed.
Probably one of the first (if not the first) periodic re-optimisations is due to \cite{psaraftis1980dynamic}.
He applied dynamic programming and periodic re-optimisations to solve a dial and ride problem.

The dial-a-ride problem (DARP) \citep{lois2017online, kirchler2013granular} is closely related to the mTSP.
It differs from the D-mTSP by having a pick-up and drop-off location.
However, those locations could be aggregated into a single abstract node; 
which allows to transform the dial-a-ride problem into an asynchronous mTSP.
Typical DARP formulations emphasis time windows.
Surprisingly, almost all DARP formulations are in the static category.
\cite{cordeau2003dial} reviewed the DARP and confirmed this view.
\cite{madsen1995heuristic} is one of the few works which
analyse the DARP in a dynamic environment intended for online scheduling.
The practical requirement was to handle up to 300 requests having 24 vehicles available to be scheduled.
Furthermore, a solution had to be returned within 2 seconds.
Multiple objectives and constraints had to be considered.
They started their work based on an algorithm by \cite{jaw1986heuristic}.
Again, their approach falls into the class of random insertion heuristics,
where a known schedule is improved.
The heuristic looks for all feasible insertions in the existing routes 
and adds the new request such that a minimal change to the objective occurs.
If a feasible insertion cannot be found the request (job) is not served.

Given the gap of a systematic dynamics' classification with the exception of information evolution and quality  
several dynamic scopes will be proposed.

\subsection{Scopes}\label{ssec:dynamics-scope}
Measures for the dynamic scopes will be introduced,
which include absolute, relative and vehicle ($m$) dependent values.

\textit{Absolute dynamics} $D_a \in \sett{n-1}$ is defined as the fixed number of 
customer requests known at time $t$ with the exception of the final period.
If $t \in T$, where $T$ is a set of equidistant time steps 
and $D_a$ is the same for each $t \in T$,
then I will call this \textit{absolute sequential dynamic}.
To accommodate the number of vehicles the term $m$-absolute dynamics is proposed: $D_a^m = \frac{D_a}{m}$.
I will also refer to this as \textit{vehicle dynamics}.
\begin{example}[Sequential absolute dynamics]
	Customers $C=\seq{1,2,\dots,100}$ ordered by ``reveal'' time $T=\seq{5,10,\dots,500}$ 
	and let the absolute dynamics be $D_a=5$.
	That means at $t_1$ five customers $\set{c_1,\dots,c_5}$ 
	and related deterministic or stochastic quality information such as their locations are known.
	So, at ``any'' point of time $t_k = 5k$ the customers $\set{c_k,\dots,c_{k+5-1}}$
	are known and the previously revealed customers  $\set{c_1,\dots,c_{k-1}}$.
	Obviously, when reaching the end ($n-D_a$) of the scenario 
	only the remaining customer information is revealed.
	
	Assume there are three vehicles and the $m$-absolute dynamics is $D_a^m=D_a^3=2$ customers. 
	This is equivalent to absolute dynamics $D_a=6$.
	Let us consider the case of sequential dynamics.
	At time $t_1=5$ customers $\set{c_1,\dots,c_6}$ are visible.
	Assume an algorithm assigned $\set{c_2,c_3,c_6}$ to the vehicles before $t_2 = 10$,
	and that these customers can be serviced within $\Delta t = 5$.
	That means at $t_2$ another three customers need to be made visible to fulfil $D_a=6$.
	Hence, at $t_2$ the customers  $\set{c_1,c_4,c_5,c_7,c_8,c_9}$ are
	available for the algorithm.
	Note that all vehicles can be used again.
	Hence, at each time step three new customers can be serviced.
\end{example}
An interesting special case occurs when $m=D_a$ and the next customer locations are randomly distributed in the ``vicinity'' 
of the last visited customers, because this is equivalent to random walks in higher dimensions.
It would be interesting to investigate the relationship between $m$ dimensional Gaussian random walks and the D$_a$-mTSP 
with the closest vehicle heuristic. This could lead to interesting insight 
and relationships to the Black-Scholes derived approaches.
Later, I will introduce an algorithm 
which solves the case $m=D_a$ optimally 
under the constraints of sequential dynamics, equidistant time steps ($\Delta t$) and $m$ customers being serviced within $\Delta t$.

\begin{example}[Ridesharing and ridehailing services]
	Uber, Lyft, Bolt and many more app-taxis (ride sharing companies) 
	have become common and are a serious competition to taxicabs (ridehailing services).
	However, they experience the same operational characteristics.
	Customers waiting for a service occurs when $D_a > m$ at time $t$.
	Customers may still wait for a service, when $D_a\leq m$ at time $t$,
	even more vehicles (servers) are available than customers.
	This is due to the random arrival of customers.
	
	In this example we are touching on the relationship between dynamic routing and queueing systems.
	In a nutshell it is possible to use queueing systems
	to explain operational characteristics such as
	customers waiting for service, average waiting time and others.
	
	Let us consider the $M/M/m$ queueing system.
	Here, customers arrive according to a Poison process (first $M$ reflects the Markov process - Poison process is a subset); 
	and are serviced within exponentially distributed (second $M$) times 
	by $m$ servers (vehicles).
	The birth-death process describes this system.
	
	We can define the arrival rate $\lambda$ as the number of customers per hour being visible to the system.
	This is equivalent to the absolute dynamics scope.
	For instance, three customers arrive per hour means $\lambda = 3\text{C/h} = D_a$.
	The (single) service rate $\mu$ is the numbers of customers one server (vehicle) can visit per hour. 
	That means, the service rate of the fleet is $m\mu$.
	Queueing theory is mainly concerned with stable systems.
	That means were the arrival rate is smaller than the service rate.
	$\rho = \frac{\lambda}{m \mu}$ is defined as the utilisation.
	So, for a stable system $\rho <1$ is required 
	otherwise the system is unstable and the queue will grow without bound assuming an infinite time horizon 
	and customers not leaving the queue before they are serviced.
	
	The subsequent test instances and continuous approximation model 
	will consider the special case of constant service and constant arrival rate with finite time horizon.
	Furthermore, we will require that $\lambda = D_a$ (except at the end of the time span).
	The (single) service rate is $\mu = D_a^m$ and the fleet service rate is $m \mu = D_a$.
	
	This example should give an idea about the large scope dynamic routing can assume.
\end{example}

If the total number of customer requests $n$ is known
then \textit{relative dynamics} can be defined:
\lequ{rel-dynamics}{D_r = \frac{k}{n},}
where $k \in \setn$ controls the fraction of revealed data.
The use of these rational numbers has the advantage that
one can convert between relative and absolute dynamics without rounding.
Alternatively, a given percentage value $q$ is converted to the relative (node) dynamics via $D_r = \frac{\round{q n}}{n} $.
Generally, it is more convenient to allow percentages as relative dynamics.
I will use \textit{customer dynamics} synonymously with relative dynamics. 
\begin{example}[Relative dynamics]
	(a) Customers $C=\seq{1,2,\dots,100}$ ordered by ``reveal'' time $T=\seq{5,10,\dots,500}$ 
	and let the relative dynamics be $D_r=\frac{k}{n}=\frac{5}{100}=.05$.
	That means 5\% of the customer data is revealed at each time step $t_k,~k\leq 96$.
	At time step 20 all the customer data is known.
	
	(b) $D_r=100\%$ is the special case of all data being known at the beginning.
	This means the BD-mTSP degenerates to the B-mTSP.
	
	(c) Assume that the target is to server 100 customers during a working day (having 10 hours).
	Furthermore, customer updates are processed every 30 minutes.
	This gives the relative dynamic scope $D_r=5\%$.
\end{example}	
To consider the number of vehicles, $m$ needs to be incorporated.
Equally to before, $m$-relative dynamic is defined:
\lequ{m-relative}{D_r^m = \frac{D_r}{m}.}
\begin{example}[$m$-relative dynamics]
	The number of vehicles is $m=3$, customers $C=\seq{1,2,\dots,100}$ are ordered by ``reveal'' time $T=\seq{5,10,\dots,500}$ 
	and let the $m$-relative dynamics be $D_r=mD_r^m=2\%$.
	That means $D_r=6\%$ of the customer data is revealed at each time step.
\end{example}	
In most cases it is necessary to convert to absolute node dynamics:
$D_a = \round{m D_a^m}$, 
$D_a = \round{n D_r}$ and $D_a = \round{n m D_r^m}$.
Generally, when vehicles are the driver to describe the dynamics scope, then $m$-absolute dynamics are used. 
For instance, $D_a^m=3$ customers per hour per vehicle with four vehicles requires a $D_a=12$ C/h scope.
When customers are the focus to describe the dynamics scope, then relative dynamics are used.
For instance, let $D_r=5\%$ with 300 customers per hour
and four vehicles operating then the absolute dynamic scope is $D_a=60$ C/h and the vehicle-dynamics are $D_a^m=15$ C/h.

The next level of dynamics occurs, when the number of customers varies at each time step $t_k$.
I will call this \textit{variable dynamics} $D_v$. 
$D_v = \seq{d_1,\dots,d_n}$ is a sequence of customers visible at $\seq{t_1,\dots,t_n}$.
Note, that at $d_n=1$ at $t_n$. 
Generally, due to the finite nature of the formulation the following set of inequalities hold:
$\set{d_1\leq n, d_2\leq n-1, \dots,d_{n-1}\leq2,d_n=1}$.
\begin{example}[variable dynamics]
	The number of vehicles is $m=3$, customers $C=\seq{1,2,\dots,100}$ are ordered by ``reveal'' time $T=\seq{5,10,\dots,500}$ 
	and let the variable dynamics be $D_v=\seq{3,4,5,3,4,5,\dots,3,4,5}$.
	That means at $t_1=5$ customers $c_1,c_2,c_3$ are visible.
	At $t_2=10$ information about customers $c_2,c_3,c_4,c_5$ is available.
	At $t_3=15$ information about customers $c_3,\dots,c_7$ is available.
	Now the visibility reduces to customers $c_4,c_5,c_6$  at $t_4=20$.
	This is interesting because it offers considerations such as allowing memory of $c_7$.
	Alternatively, once visibility is lost it could be that the location of $c_7$ changed.
\end{example}

As mentioned before the most frequently used type of dynamic is random insertion.
In essence $\seq{c_1,\dots,c_k}$ customers are known and $c_{k+1}$ is added.
Hence, this is a special case of the absolute dynamic scope with $D_a\geq1$.
It should be noted that in the previously reviewed literature time windows for customer tasks were considered.
This adds another timeline to the problem, or adds constraints to the existing timeline $T$.
To distinguish these two the term \textit{event-knowledge-timeline} and \textit{scheduled-timeline} are proposed.
Future work shall look at the dynamics of the mTSP with time-windows.

On top of these fundamental dynamics the stochastic elements of customers' quality 
of the location (e.g. position improves over time) must be taken into account.
There are two aspects - where and when. 
A trivial case is that the customer location is static over the entire period.
The other extreme is that each customer's location changes continuously over time (``moving target'').
A potential approach is allowing on top of inserting customers the removal of customers.
Additionally, the accuracy of the static or changing location of each customer needs to be considered.
Now, these measures could apply to classes of customers or the entire customer set.

Vehicle dynamics can be considered using a varying number of active vehicles over time.
Additionally, uncertainty in the position of the vehicles position or state can occur as well.

\begin{figure}[htbp]\center
	\includegraphics[width=\columnwidth, height=6cm, keepaspectratio]{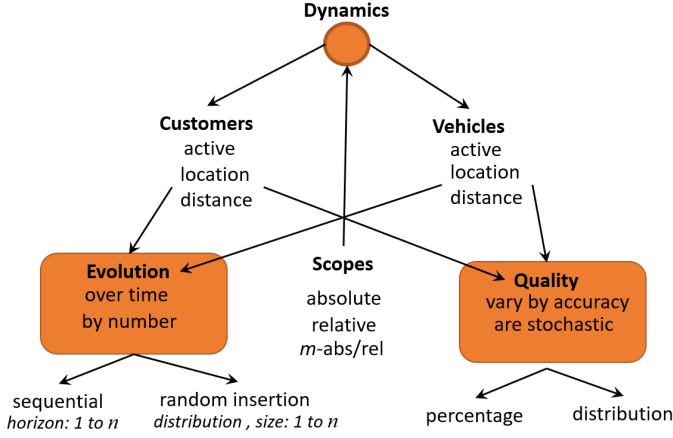}%
	\caption{Types, factors, and scopes of dynamics in mTSP.}%
	\label{fig:dynamics}%
\end{figure}

Figure \ref{fig:dynamics} gives an overview of the above-mentioned dynamics that can occur in the mTSP.

%% file: sec-Heuristics.tex
\section{BD-mTSP Heuristics}\label{sec:Heuristics}
We will consider two greedy balanced-dynamic mTSP (BD-mTSP) heuristics. 
The first heuristic allocates the vehicle closest to ``visible'' customers. 
The second heuristic assigns vehicles minimising their distance.

The input for the algorithms is a 
distance matrix $D$, 
the number of vehicle $m$
and the absolute dynamics $d$.
The balancing threshold (capacity limit) is determined automatically.
Please note that the algorithms can be easily amended to contain the capacity as input and $m$ as variable.
The output of the algorithms is a list of $m$-routes.

The dynamic algorithms are examined using and introducing several test instances.
The difference of the dynamic results to the static exact mTSP test instances are explained.
Later (Section \ref{sec:CAM}), the heuristics will be used to derive continuous approximation models (CAMs).

\subsection{Closest Vehicle\label{ssec:BD-CVH}}
This algorithm is the dynamic version of the CVH.
I introduced the CVH in \cite{garn2020closed}.
This algorithm uses absolute customer dynamics,
i.e. $d = D_a$ is the number of customers visible at $t_k$.
As shown in the previous section relative dynamics and $m$-dynamics measures can be converted to absolute dynamics.
It should be noted that the algorithm can be easily adapted to variable customer dynamics 
by providing a vector $d$ as input instead of a single value (and adapting Line 6).

$m$ vehicles are available with a capacity of $Q=\ceil{\frac{n-1}{m}}$. 
The customer nodes are balanced between those vehicles using $Q$.
This is a strict balancing threshold.
In case $\mod (n-1,m)=0$ the vehicles are perfectly balanced,
otherwise $\mod (n-1,m)=r$ customers are serviced by one or more vehicles.
It may be preferable to allow $Q$ as input to the algorithm to allow more flexibility.
Providing a capacity limit for each vehicle is the most versatile variant.

We will assume that $d$ customers and their locations are revealed sequentially,
i.e. at each time step $d$ customers are visible.
Furthermore, we will require $\min \set{m,d}$ vehicles to have completed their assigned tasks (customer visits).
This constraint ensures that all vehicles are active at each time step,
otherwise it could happen that one visits $Q$ customers in sequence, 
before the next vehicle starts. 
Hence, this means the number of customers serviced by each vehicle is balanced at each time step 
(unless $D_a<1$).

It is assumed that the total number of customers $n$ is known in advance, 
but their locations are not.
As an option, $n$ could be estimated using $Q=\ceil{\frac{n-1}{m}}$, 
i.e. $n=mL+1$ if $\mod(n-1,m)=0$ otherwise $n=Qm-m+r+1$ (with $\mod(n-1,m)=r$).
At each time step, we will require that $\min \set{d, m}$ customers 
are serviced, and that a vehicle cannot service more than one customer.

Algorithm \ref{alg:BD-CVH} shows the details for this implementation.
\begin{algorithm}[!h]
	\caption{Dynamic Closest Vehicle Heuristic. \label{alg:BD-CVH}}
	\begin{algorithmic}[1]
		\Require distance matrix $D=(d_{ij}) \in \real^{n \times n}$; number of vehicles $m \in \setN$; dynamics $d \in \setN$
		\Ensure routes $r = \seq{r_1,\dots,r_m}$
		\State $f=e_m=\mat{1&1&\dots&1}$ \Comment{all vehicle start from 1} 
		\State $r=e_m$ \Comment{initialise routes}
		\State $v=\set{2,\dots,n}$ \Comment{nodes not visited}
		\State $y = \mat{0,0,\dots,0}$ \Comment{number of nodes visited per vehicle}
		\While { $\#v > 0$}
			\State $\delta = \set{v_k ~|~ k \leq d}$ \Comment{visible nodes}
			\State $\Delta = D_{f,\delta}$ \Comment{sub-matrix}
			\For {$e \in \set{1,\dots,\min \set{m,|\delta|}}$}
				\State $\mat{i&j} = \argmin \Delta$; \Comment{closest vehicle}
				\State $x_{ij}=1$; \Comment{assigned vehicle}
				\State $\Delta_{i:}=\infty,~\Delta_{:j}=\infty$ \Comment{inaccessible row \& column}
			\EndFor
			\State $t = X \delta$; \Comment{to-nodes}
			\For {$k \in m$} ~\Comment{fore each vehicle}
				\If {$t_k>0$} ~\Comment{vehicle active}
					\State $r_k = r_k \cup t_e$ \Comment{add to route}
					\State $y_k = y_k + 1$ \Comment{increment visited nodes}
					\If{$y_k \geq \ceil{\frac{n-1}{m}}$} \Comment{balancing condition}
						\State $D_{t_k:} = \infty$ \Comment{vehicle reached limit}
					\EndIf
				\EndIf
			\EndFor
			\State $v = v \setminus t$ \Comment{remove used to-nodes}
			\State $f_k = t$; \Comment{new from nodes}
		\EndWhile
	\end{algorithmic}
\end{algorithm}
It begins with requiring some input such as $D$.
In real-world implementations $D=(d_{ij})$ would not be known in advance.
For road networks shortest path algorithms can be used to derive $D$.
In most of the subsequent test instances customer locations were given 
and the Euclidean distance between customer location $i$ and $j$ are derived using $d_{ij}=\sqrt{(x_i-x_j)^2+(y_i-y_j)^2}$.
Sometimes (e.g. for extra-large test instances) it may be more practical or necessary to compute $D$ during execution.
This can be achieved in Line 7.

The algorithm returns routes for all vehicles.
Again, in an online version it is recommended to use the route output from Line 16.

The first five lines of the algorithm are initialisation and an overall loop starts.
Line 6 controls the nodes visibility (information evolution).
This represents the limitations imposed by sequential absolute dynamics.
Consequently, the distance matrix under consideration is reduced to $\Delta$ (Line 7).
As mentioned before, for real-world implementations, 
this is the point where the new customer locations are known (rather than revealed) and $\Delta$ 
is computed instead of being extracted from $D$.
Although pre-processed distance matrix constitutes computational savings.
However, for the study of dynamics this implementation detail does not make a difference.
Line 8-12 assign the closest vehicle to the visible nodes.
Line 9 returns the row $i$ and column $j$ of the minimum matrix element.
This is equivalent of identifying the closest vehicle.
In case of identical distances the first found index is returned.
In Line 10 that vehicle is assigned to the customer ($x_{ij}=1$) using the assignment matrix $X \in \setB^{m\times |\delta|}$.
Note, it can happen that $|\delta|<m$.
In order to prevent the choice of the same vehicle and customer the respective row and column are set to $\infty$.
Line 13 returns the original node numbers by making use of the matrix multiplication.
Line 14-22 add the vehicle to the route and enforce the balancing constraint.
Line 23 removes the visited nodes from all not visited nodes.
Syntactically, the case $t_k=0$ may need to be considered (depending on the programming language).
Line 24 updates the from-nodes with the to-nodes used in the last step, 
leaving the unused vehicles (some from nodes) as they are.

\subsection{Dynamic Assignment Vehicle Heuristic (BD-AVH)\label{ssec:BD-AVH}}
This algorithm builds on the conceptual framework of the BD-CVH.
Line 8-12 were responsible for assigning the closest vehicle,
which is a greedy heuristic. 
A better solution is obtained by assigning available vehicles such 
that the distance is minimised.
This is achieved with the following binary program (if $m\leq d$):
\LP	{\label{bp:a1}\obj[min]{\sum_{i=1}^m \sum_{j=1}^d \Delta_{ij} x_{ij} \mcomment{assignment distance}}}
	{\sto{\sum_{j=1}^m x_{ij} = 1,~i \in \sett{d} \mcomment{customers}}
	 \sto{\sum_{i=1}^n x_{ij} \leq 1, ~j\in \set{1,\dots,m} \mcomment{vehicles}}
	}
In the case of $m>d$ the constraints change slightly, and the following program is solved instead.
\LP	{\label{bp:a2}\obj[min]{\sum_{i=1}^m \sum_{j=1}^d \Delta_{ij} x_{ij} \mcomment{assignment distance}}}
{\sto{\sum_{j=1}^m x_{ij} \leq 1,~i \in \sett{d}  \mcomment{customers}}
	\sto{\sum_{i=1}^n x_{ij} = 1, ~j\in \set{1,\dots,m} \mcomment{vehicles} }
}
Consequently, the algorithm for the BD-AVH is the same as the BD-CVH apart from line 8-12 being replaced with the above binary programs, 
i.e.  if $m\leq d$ then program (\ref{bp:a1}) else program (\ref{bp:a2}).
I will call this part of the algorithm \textit{vehicle assignment}.

\begin{proposition}(D-AVH optimal algorithm for $d \leq m$)
The D-AVH is an optimal algorithm when $d$ vehicles are required to be assigned 
if and only if the dynamics are sequential and absolute, $d \leq m$, and customer locations beyond the dynamics scope $m$ are unpredictable.
\end{proposition}
If $d \leq m$ then the vehicle assignment part of Algorithm \ref{alg:BD-CVH} uses the binary program \ref{bp:a1},
which is optimal. 
Since, the dynamic scope is sequential and locations beyond the dynamics scope $m$ are unpredictable this is
the best solution that can be obtained.
Hence, there is no algorithm that is better. 
Note, that other algorithms can find better final route solution values; but this is coincidental.
Since the vehicle assignment phase is optimal no other algorithm can find a better solution value,
because we assumed that the customer locations beyond scope $m$ are unpredictable.

Note that when the requirement of assigning $d$ vehicles is omitted 
then it is not guaranteed that D-AVH is optimal. 
For instance, it may be better if a single vehicle visits the $d$ customers.
However, assigning $d$ vehicles has the advantage that they are balanced, 
and $d$ customers are serviced at each step. 
It is also interesting to note that when $d>m$ a potentially better solution may be found
by building partial anticipating paths.

Implementing meta heuristics instead of the vehicle assignment can improve 
the solution quality when $D_a > 1$. 
One of the most popular heuristics for the CVRP is the Iterated Local Search (ILS) one.
This meta-heuristic was successfully applied by \cite{Penna2013} and \cite{uchoa2017new}.
It seems to outperform most other heuristics in terms of solution quality and computational time. 
Another interesting ILS implementation was presented in the context of cross-docking and the VRP by \cite{morais2014iterated}.
Hence, adapting this heuristic and replacing the vehicle assignment program would enhance the solution quality.
Other meta heuristics such as greedy randomised adaptive search procedure (GRASP),
Tabu Search, and
Biased random-key genetic algorithms (BRKGA) are also likely candidates to offer improved solution values.
BRKGA was used by \cite{Ruiz2019} to solve open vehicle routing problems, and introduced in a more general context by \cite{Goncalves2011}.
\cite{kulak2012joint} used Tabu Search with a clustering algorithm in a related context.
Their work will be considered in Section \ref{sec:real-world-apps}.

The next subsection gives insights about the behaviour of the dynamic scopes 
and the performance of the heuristics.
\input{ssec-TestInstances}

%% file: ssec-TestInstances.tex
\subsection{Test Instances} \label{ssec:heuristics-test-instances}
The balanced dynamic closest vehicle heuristic (BD-CVH) and
 balanced dynamic assignment vehicle heuristic (BD-AVH) are examined
using three sets of test instances. 
A fourth set with extra-large test instances highlights computational challenges.

The first set introduces new test instances, and 
some were adapted from TSPLIB (\url{http://elib.zib.de/pub/mp-testdata/tsp/tsplib/tsp/}).
This set considers micro to medium customer range test instances.
Table \ref{tab:diff-CVH-AVH} shows the class name and associated customer range.
The number of vehicles $m$ for this set was chosen arbitrarily varying between 2 and 20.

The second set of mTSP test instances were created using the well-known E and K set often used for TSP analyses.
These test instances were adapted from the TSP instances Christofides/Eilon \citep{christofides1972algorithms} and Krolak/Felts/Nelson by adding multiple vehicles and dynamics.
The number of customers in this set varies between 51 and 200. 
That means, the set contains small to medium test instances.
The number of vehicles $m$ were set to be between 2 and 5.

The third set of test instance originated from the Capacitated Vehicle Routing Problem Library (CVRPLIB, \url{http://vrp.atd-lab.inf.puc-rio.br}).
It contains Uchoa CVRP instances \citep{uchoa2017new} 
with customers varying between 115 and 936 where demand is unitary.
That means, the set contains medium to large test instances.

The three sets and more are available on \url{http://www.smartana.org}.

$m$-absolute sequential dynamics (vehicle dynamics) and 
relative sequential dynamics (customer dynamics) were considered.
The set $\set{\frac{1}{2},1,\frac{3m}{2},2,4,8}$  is used for the $m$-absolute sequential dynamics.
$D_a^m=\frac{1}{2}$ means that half of the fleet is not assigned at each time step.  
When there is exactly one customer available for a vehicle then $D_a^m=1$. 
In case there is more than one customer available for each vehicle at each time step $D_a>1$.

The chosen relative dynamics $D_r$ are 2\%, 5\%, 7\%, 10\%, 20\%, 30\% and 100\%.
Relative dynamics are based on the total number of customers $n$.
The absolute dynamics can be obtained by $\round{nD_r}$.

Table \ref{tab:set1-rel} shows the results for the first set of test instances.
\input{tab-set1-rel}
The instance \code{garn9-m2} shows that both BD-CVH and BD-AVH are identical.
It is also interesting to compare them to their static versions,
which I proposed in \cite{garn2020closed}. 
I used the same test instances and the results are shown in \citep[p4, Table 1]{garn2020closed}.
Note that when dynamics $D_r=100\%$ the BD-CVH must give the same results as the static B-CVH.
The optimal static mTSP solution value is 44.8, which the dynamic version achieves with $d=1$.
The reason for this is that the order of the provided customers matters.
For the micro test instances ($n \in (0,30]$) some of the dynamics are too small 
to generate solutions, because $\round{nD_r} = 0$.
Trivially, the dynamic solution values for the \code{garn13-m3L4} and \code{garn20m3} instances 
are higher than the optimal solution values of the static mTSP which are 57.7 and 90.5.

Generally, we would expect higher node visibility to lead to better results when $D_a \geq  m$.
This theory will be confirmed with the continuous approximation model (Section \ref{sec:CAM}).
On rare occasions the BD-CVH returns better solutions, e.g. instance \code{eucl-n100m7} with $D_a=20$.
This is justifiable by observing that although the assignment step
is optimal the overall algorithm is a heuristic. 
Hence, the assignment can send the vehicles to an unfavourable location for subsequent customer destinations.

Lower solution values for instances with $D_a<m$ than the $D_a=m$ instance are expected, 
because of the increased flexibility.
This can be observed in Table \ref{tab:set1-abs}, \ref{tab:set2-abs} and \ref{tab:setX-abs}. 
The highest solution values are at $D_a^m=1$.
The phenomenon that $D_a^m=1$ has the highest solution value 
is due to the reduced flexibility.
\begin{proposition}[$D_a^m=1$ maximum solution value]
	The dynamic mTSP solution value is at a maximum when $D_a^m=1$, 
	given absolute sequential dynamics.
	That means, $D_a^m=1$ implies a solution value $s^* > s_m,~ s_m \in S$,
	where $S$ is the set of solutions obtained with an optimal dynamic mTSP algorithm.
\end{proposition}
When $D_a^m<1$ then there is a choice of vehicles to be assigned to the customers at each time step.
Hence, an optimal dynamic mTSP algorithm must return the same solution value or a better one.
Assume that $D_a^m>1$ then there is a choice of customers to be allocated to vehicles.
Hence, the increased choice must lead to solution values which are greater or equal to $s^*$.
However, the highest solution value does not necessarily imply that $D_a^m=1$, 
because $D_a^m\neq1$ can lead to similar solution values.

In general, it can be observed that the BD-AVH delivers shorter routes than the BD-CVH. 
Table \ref{tab:diff-CVH-AVH} provides details of the differences for
micro, small, medium, and large test instances.
\input{tab-diff-CVH-AVH}
The three sets of test instances were used for $m$-absolute and relative dynamics.
Let $L_i^a$ and $L_i^c$ be the total distance for test instance $i$ obtained from the BD-AVH and BD-CVH routes, respectively.
For each test instance the relative difference $\delta_i$ was computed by
\lequ{rdti}{\delta_i = \frac{L_i^a - L_i^c}{L_i^a}.}
For instance, 42 large $m$-absolute dynamics test instances ($n\in(400,1000]$) were used
to determine relative differences ($\delta_i,~i \in \sett{42}$).
The average relative difference is 5.0\%, i.e. BD-CVH has five percent longer routes.
On average, across all test instances the BD-AVH has shorter (1.8\% to 5.0\%) total route lengths than the BD-CVH.
However, this comes with the cost of 97\% longer run-times.
The time-percentage is determined similar to Equation \ref{eq:rdti}.
The average run-times for the BD-AVH are 52.8, 193.1 and 365.3 milliseconds for instances in set 1, 2 and 3.
This was measured on a computer with an 8 Intel i9@3.6GHz processors and 32GB RAM.
Hence, the algorithms are suitable for online processing.

A fourth set of extra-large test instances with customers beyond the 1000 threshold was considered.
However, due to computational challenges not sufficient results were obtained.
I will briefly discuss the proposed set and the challenges.
The set originated from \cite{arnold2019efficiently}.
In the original test instances, the number of customers $n$ vary between 3,000 and 30,000; 
demand $q$ varied; and each fleet was homogeneous with a capacity limit $q^*$.
These instances were transformed to be suitable for the balanced mTSP by:
(a) setting the demand to a unitary one;
(b) computing the number of vehicles: $m = \floor{1.05 \frac{\sum_{i=1}^n q_i}{q^*}}$; and
(c) readjusting the load's upper limit of vehicles to $Q=\ceil{\frac{n}{m}}$.
The first challenge occurs when determining the distance matrix.
The number of customers $n=30k$ leads to a matrix with $n=900M$ entries, 
which uses $8n$ bytes. 
That means, at least 7.2GB are required for a single test instance of size $n=30k$.
This issue can be ``partially'' resolved by adapting the algorithm.
Instead of requiring $D$ as input the coordinates are provided.
Note, once the dynamic TSP becomes static or reaches a memory threshold, further adaptations need to be done 
(e.g. breaking up the distance matrix).
The second issue happens when using the Assignment vehicle heuristic.
Computing the optimal assignment requires a matrix of constraints of size $m \times m d$.
For instance, $m=359$ and $d=180$ is the dimension of one of the first assignment problems to be solved.
That means, the problem has $m\times d=64,620$ decision variables.
This leads to a matrix of 23.2 million entries necessary for the customer constraints.
Hence, the assignment optimisation problem for extra-large instances cannot be solved using standard optimisation software.
As a consequence, the BD-CVH was used instead of the BD-AVH.
Table \ref{tab:setXXL-abs} shows the results for the solved extra-large test instances.
Generally, the results follow the insights suggested by this paper's propositions 
even though the BD-CVH algorithm was used instead of the BD-AVH.
The runtime is surprisingly fast between 0.6 seconds ($n=11k$, $m=115$, $D_a^m=\frac{1}{2}$) and 7.65 minutes ($n=20k$, $m=717$, $D_a^m=8$).

In summary this section introduced two heuristics, 
which can be used for online routing.
Three sets of test instances were analysed.
In general, the BD-AVH returns better solution values than the BD-CVH but has longer run-times.
D-AVH is an optimal algorithm for $m\leq d$.
Computational challenges when using the D-AVH for extra-large test instances were explained, 
and BD-CVH solutions were provided.

%% file: tab-set1-rel.tex
\begin{table*}[htbp]
	\centering
	\caption{BD-CVH and BD-AVH test instances for set 1 using relative dynamics.}
	\begin{tabular}{c|ccc|ccc|c}
		\toprule
		\boldmath{}\textbf{$D_r$}\unboldmath{} & \textbf{2\%} & \textbf{5\%} & \textbf{7\%} & \textbf{10\%} & \textbf{20\%} & \textbf{30\%} & \textbf{100\%} \\
		\midrule
		garn9-m2 & \multicolumn{2}{c|}{$D_a = 0$} & \multicolumn{2}{c|}{$D_a = 1$} & \multicolumn{2}{c|}{$D_a = 2$} & $D_a = 8$ \\
		\midrule
		BD-AVH &      & \multicolumn{1}{c|}{} & \multicolumn{2}{c|}{44.8} & \multicolumn{2}{c|}{68.9} & 44.8 \\
		BD-CVH &      & \multicolumn{1}{c|}{} & \multicolumn{2}{c|}{44.8} & \multicolumn{2}{c|}{68.9} & 44.8 \\
		\midrule
		garn13-m3L4 & \multicolumn{1}{c|}{$D_a = 0$} & \multicolumn{3}{c|}{$D_a = 1$} & $D_a = 2$ & $D_a = 4$ & $D_a = 12$ \\
		\midrule
		BD-AVH & \multicolumn{1}{c|}{} & \multicolumn{3}{c|}{71.1} & 69.2 & 91.9 & 79.4 \\
		BD-CVH & \multicolumn{1}{c|}{} & \multicolumn{3}{c|}{71.1} & 69.2 & 92.9 & 87.8 \\
		\midrule
		garn20-m3 & \multicolumn{1}{c|}{$D_a = 0$} & \multicolumn{2}{c|}{$D_a = 1$} & $D_a = 2$ & $D_a = 4$ & $D_a = 6$ & $D_a = 19$ \\
		\midrule
		BD-AVH & \multicolumn{1}{c|}{} & \multicolumn{2}{c|}{103.1} & 108.1 & 117.1 & 119.3 & 116.3 \\
		BD-CVH & \multicolumn{1}{c|}{} & \multicolumn{2}{c|}{103.1} & 111.5 & 120  & 120.9 & 120.2 \\
		\midrule
		bays29-m4 & \multicolumn{2}{c|}{$D_a = 1$} & $D_a = 2$ & $D_a = 3$ & $D_a = 6$ & $D_a = 8$ & $D_a = 28$ \\
		\midrule
		BD-AVH & \multicolumn{2}{c|}{5315} & 4114 & 4319 & 4131 & 3983 & 3458 \\
		BD-CVH & \multicolumn{2}{c|}{5315} & 4148 & 4665 & 4260 & 4136 & 3643 \\
		\midrule
		berlin52-m5 & $D_a = 1$ & $D_a = 3$ & $D_a = 4$ & $D_a = 5$ & $D_a = 10$ & $D_a = 15$ & $D_a = 51$ \\
		\midrule
		BD-AVH & 16.4k & 22.6k & 22.3k & 25.7k & 17.1k & 15.2k & 13.6k \\
		BD-CVH & 16.4k & 23.2k & 22.6k & 26.8k & 18.2k & 15.9k & 13.6k \\
		\midrule
		eucl-n100m7 & $D_a = 2$ & $D_a = 5$ & $D_a = 7$ & $D_a = 10$ & $D_a = 20$ & $D_a = 30$ & $D_a = 99$ \\
		\midrule
		BD-AVH & 39.4k & 37.1k & 42.8k & 39.3k & 30.6k & 29.9k & 25.8k \\
		BD-CVH & 38.9k & 39.5k & 44.7k & 38.4k & 29.8k & 30.9k & 26.4k \\
		\midrule
		lin318-m20 & $D_a = 6$ & $D_a = 16$ & $D_a = 22$ & $D_a = 32$ & $D_a = 63$ & $D_a = 95$ & $D_a = 317$ \\
		\midrule
		BD-AVH & 222.2k & 297.5k & 341.1k & 319.6k & 252.3k & 202.6k & 199.1k \\
		BD-CVH & 224.9k & 306.3k & 349.1k & 331.1k & 263.3k & 207.8k & 202.6k \\
		\bottomrule
	\end{tabular}%
	\label{tab:set1-rel}%
\end{table*}%

%% file: tab-diff-CVH-AVH.tex
\begin{table}[htbp]
	\centering
	\caption{Difference between BD-CVH and BD-AVH.}
	\adjustbox{width=\columnwidth}{ 
	\begin{tabular}{cc|ccc|ccc}
		\toprule
		\multicolumn{2}{c|}{\textbf{test instance}} & \multicolumn{3}{c|}{\textbf{$m$-absolute dynamics}} & \multicolumn{3}{c}{\textbf{relative dynamics}} \\
		class & \#customers & mean &  count  & median & mean &  count  & median \\
		\midrule
		micro & (0, 30] & 2.6\% & 24   & 1.2\% & 2.4\% & 18   & 1.2\% \\
		small & (30, 100] & 1.8\% & 84   & 1.7\% & 2.7\% & 98   & 2.3\% \\
		medium & (100, 400] & 3.3\% & 132  & 2.6\% & 2.5\% & 154  & 1.9\% \\
		large & (400, 1000] & 5.0\% & 42   & 4.9\% & 2.3\% & 49   & 2.1\% \\
		\bottomrule
	\end{tabular}%
	}%
	\label{tab:diff-CVH-AVH}%
\end{table}%

%% file: sec-real-world-apps.tex
	\section{Real-world applications \label{sec:real-world-apps}}
	This section shows how to apply the BD-mTSP to real world scenarios.
	Hence, it will propose steps to transform real-world scenarios 
	to balanced-dynamic multiple travelling salesmen problems.
	
	Three types of applications are considered: deliveries, warehouse transfers and taxi services.
		We begin with the well-known CVRP - Fisher's instances \cite{fisher1994optimal}, 
	which focus on deliveries.
	However, they require a bit of modification in regarding demand and dynamic.
	The second type of applications deals with transfers of pallets within warehouses,
	which motivated this research in the first place.
	The last application is about taxi services (ridehailing/sharing),
	which is relevant for providers such as Uber, Lyft, Bolt, Grab and many more.
		
	\subsection{Fisher instances}
\begin{figure*}[htbp]\center
		\includegraphics[width=\textwidth, height=6cm, keepaspectratio]{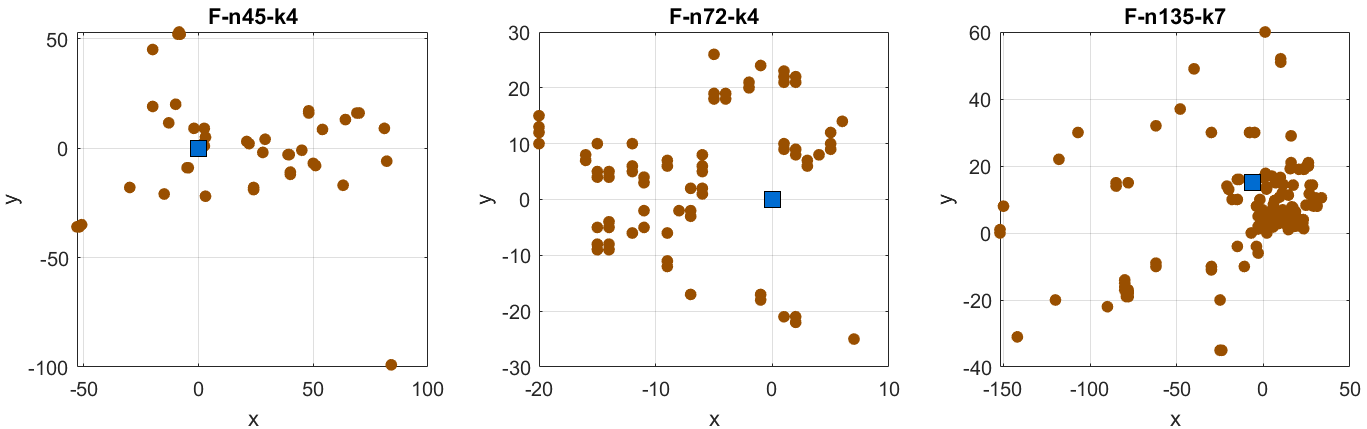}%
		\caption{Warehouse asymmetric TSP example.}%
		\label{fig:Set-F-loc}%
	\end{figure*}

	Fisher gave several real-world instances, which can be formulated as CVRP.
	Instances F-n45-k4, F-n72-k4 and F-n135-k7 are considered here (see Figure \ref{fig:Set-F-loc}).

	\textbf{Instance F-n45-k4 }describes grocery deliveries by National Grocers Limited. 
	Four vehicles ($m=4$)  begin their journey in Peterboro (Ontario terminal)
	and visit $ n-1=44 $ customers within a day.
	Customer demand and vehicle capacity exist in the data. 

	The CVRP is converted into a BD-mTSP 
	by assuming that delivering goods to the customer has priority over demand.
	This allows us to replace various demands with the requirement of visiting a customer.
	The vehicle capacity is $ L=12 $ customers.
	The Fisher instance does not specify any dynamics originally.
	The dynamics are assumed to be $D_a^m=3$ customers per vehicle - 
	reflecting schedule updates every 2 hours (assuming a working day with 8 hours). 
	The given sequence of delivery locations is permuted. 
	This permutation reflects the time-sequence of customer orders and removes unintended patterns within the original data.
	For reproducibility the original Fisher instances and the BD-mTSP dataset
	are provided on \href{https://www.smartana.org/blogs/VRP-instances.html}{smartana.org}.

	The BD-mTSP application's aim is to fulfil delivery orders as soon as possible after their occurrence. 
	Furthermore, the fleet must have its workload balanced in respect to the number of customers.
	It is assumed that each vehicle has sufficient standard groceries for its customer.
	Otherwise, intermediate vehicle stocking would be needed (or returns to the depot).
	
	The original CVRP optimal solution value is 724, when all orders are known in advance.
	The BD-mTSP solution value obtained with the BD-AVH is 1,348.
	This is 1.9 times larger than the CVRP, 
	which raises the question of how to value the increased speed in delivery versus minimising travel cost.

	\textbf{Instance F-n72-k4} is about the delivery of tires, batteries, and accessories to gasoline service stations.
	The data is associated with Exxon and formatted for CVRP.
	We assume each van has a standard stock (capacity) that can cover the needs of $L=18$ stations.
	There are $m=5$ vans and $n-1=71$ customers.
	We will consider customer visibilities between one and six customers per vehicle, i.e. $D_a^m \in \sett{6}$.
	As before a permutation was applied to the original order of the customers,
	which represent the time-sequence (or priority) of incoming orders.
	The BD-AVH obtained the following solution values:
	$L \in \set{1138, 649, 550, 441, 413, 412}$ corresponding to 
	$D_a^4 \in \sett{6}$.
	Hence, the balanced-dynamic mTSP has distances between 1.7 and 4.8 times larger than 
	its associated CVRP.
	The distances reduce with increased vehicle-dynamics (number of visible orders).
	 
	\textbf{Instance F-n135-k7} is about delivery groceries using the same company and depot as in F-n45-k4.
	This time $n-1=134$ customers are visited by $m=7$ vehicles.
	The capacity is set to $L=20$ customers per vehicle.
	The dynamics are assumed to be $D_a^m=3$ customers per vehicle.
    Again, a permutation was applied to the original order of the customers.
    The original CVRP optimal solution value is 1,162.
    The BD-mTSP solution value obtained with the BD-AVH is 2,963, which is 2.5 times larger than the CVRP.
    
    These Fisher instances gave an idea 
    about the magnitude of distance difference between static CVRPs and balanced-dynamic mTSPs.
    As a consequence, increased travel cost have to be taken into consideration, when allowing faster order fulfilment.
	 
	\subsection{Warehouses \label{ssec:wh}}
	The transfer of palettes between storage locations is an essential task in warehouses.  
	This challenge motivated this research.
	To be more precise a fleet of autonomous forklift-trucks 
	obtained transfer requests on-the-fly, which have to be fulfilled as soon as possible.
		
	Before, going into details I will mention a few related challenges.
	It is well known that order picking is one of the classic TSP applications.
	The order picking problem is part of the Steiner TSP.
	\cite{ratliff1983order} provide a solution for this problem.
	The problem we are addressing has similarities to the 
	joint order and picker routing problem \cite{kulak2012joint}.
	It can also be found within the cross-docking context.
	For instance, \cite{morais2014iterated} provided an iterated local search heuristic
	for the VRP with Cross-Docking.

\begin{figure}[htbp]\center
		\includegraphics[width=\columnwidth, height=6cm, keepaspectratio]{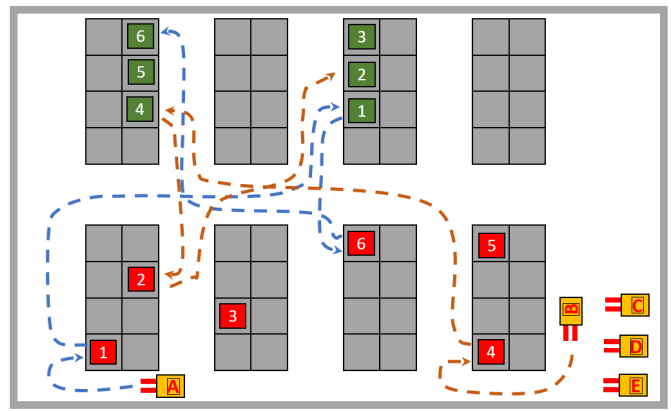}%
		\caption{Warehouse example.}%
		\label{fig:warehouse}%
	\end{figure}

	Now, let us return to the transfer challenge. 
	In a warehouse items have to be transferred between storage locations.
	Figure \ref{fig:warehouse} shows a schematic of a warehouse.
	In the example two forklift trucks are used to transfer pallets
	from the red storage locations to the green locations.
	The distance to complete the transfer jobs needs to be minimised,
	whilst ensuring that both forklift trucks are operating.
	The figure depicts a transfer job as a red and green square with identical numbers.
	For instance, the pallet at the red square labelled one needs to be
	transferred to the storage location highlighted with the green square labelled one.
	Hence, in this example there are six transfer jobs.
	Furthermore, we assume that both trucks job completion time is synchronised.
	That means, jobs are completed at the same time.
	I have started to sketch out the beginning of a feasible solution.
	In this solution forklift A travels to the storage location labelled one within a red square $ r_1 $
	and picks up the pallet from a shelf.
	The forklift travels to the storage location labelled one within a green square $ g_1 $ 
	using the shortest possible path.
	This completes one of the six transfer jobs.
	Forklift A continues from $g_1$ to $r_6$ picks up a pallet and delivers it to $g_6$.
	Forklift B starts at the same time as A.
	Its first job is $r_4 \rightarrow g_4$.
	This is followed by $r_2 \rightarrow g_2$.
	Currently, job 3 and 5 have not been allocated yet.
	
	If all jobs are known in advance, then optimal transfer sequence can be obtained using a static-mTSP.
	The key idea is to recognise that each job has a start and end location.
	Hence, a transfer job constitutes a ``logical'' node rather than a topological one.
	This permits to derive an asymmetric distance matrix as input for the mTSP.
	In the context of warehouses, distances can be computed using a shortest path algorithm.
	The transfer job itself also requires a shortest path to be solved.
	However, a transfer job's travel time (including pick-up and drop-off the pallet) does not need to be considered in the mTSP objective function.
	
	Three examples will clarify the required steps. 
	The first one illustrates the general procedure and introduces warehouse details.
	The second example extends the warehouse network 
	and visualises the entire solution route given a static asymmetric TSP.
	The third example builds on the previous two examples but uses 
	two trucks, dynamics, and balancing.
	
	\begin{example}[Asymmetric TSP in warehouse context]
		A single forklift must execute jobs 1 to 6 as depicted in Figure \ref{fig:warehouse2}.
			\begin{figure}[htbp]\center
			\includegraphics[width=\columnwidth, height=6cm, keepaspectratio]{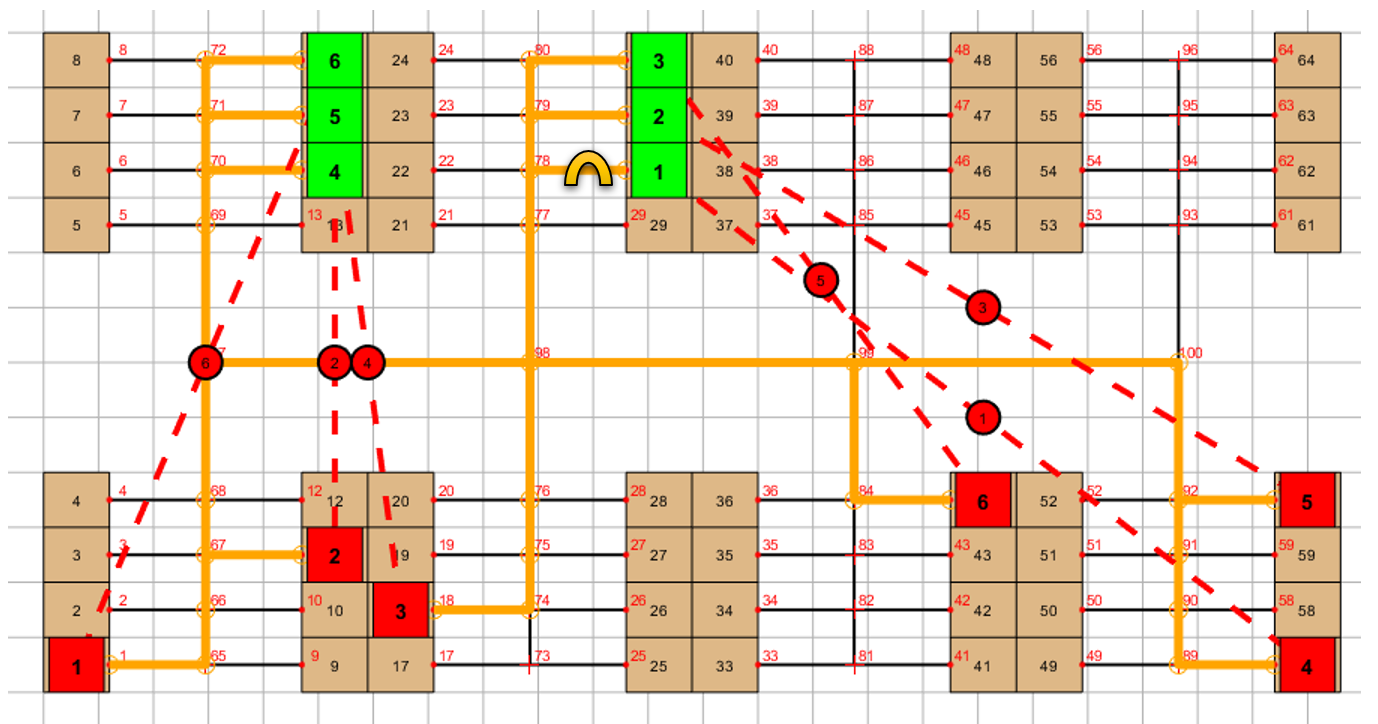}%
			\caption{Warehouse asymmetric TSP example.}%
			\label{fig:warehouse2}%
		\end{figure}
		Assume that $g_1$ is the start location (depot).
		The procedure to solve this example is:
		\begin{enumerate}
			\item Design or import warehouse layout;
			\item Create a network that connects all storage locations;
			\item Read (or provide) a job list, with pallet’s source and destination; and
				 determine shortest path for each job between its source and destination;
			\item Find shortest paths between jobs, i.e. create an asymmetric distance matrix
			\item Solve travelling salesman.
		\end{enumerate} 
	
	The wide orange lines show the topological solution, i.e. the actual travelled paths (with many overlaps).
	The dashed red lines depict the ``logical'' solution for the asymmetric TSP,
	i.e. $g_1-4-2-5-3-6-1$. 
	\end{example}
	
	\begin{example}[Exact solution of asymmetric TSP in warehouse]
		Figure \ref{fig:warehouse-depot}a shows a warehouse with 301 storage locations,
		ten transfer-jobs and one forklift truck located in the depot. 
		\begin{figure*}[htbp]\center
			\begin{tabular}{cc}
				\includegraphics[width=.5\textwidth, height=9cm, keepaspectratio]{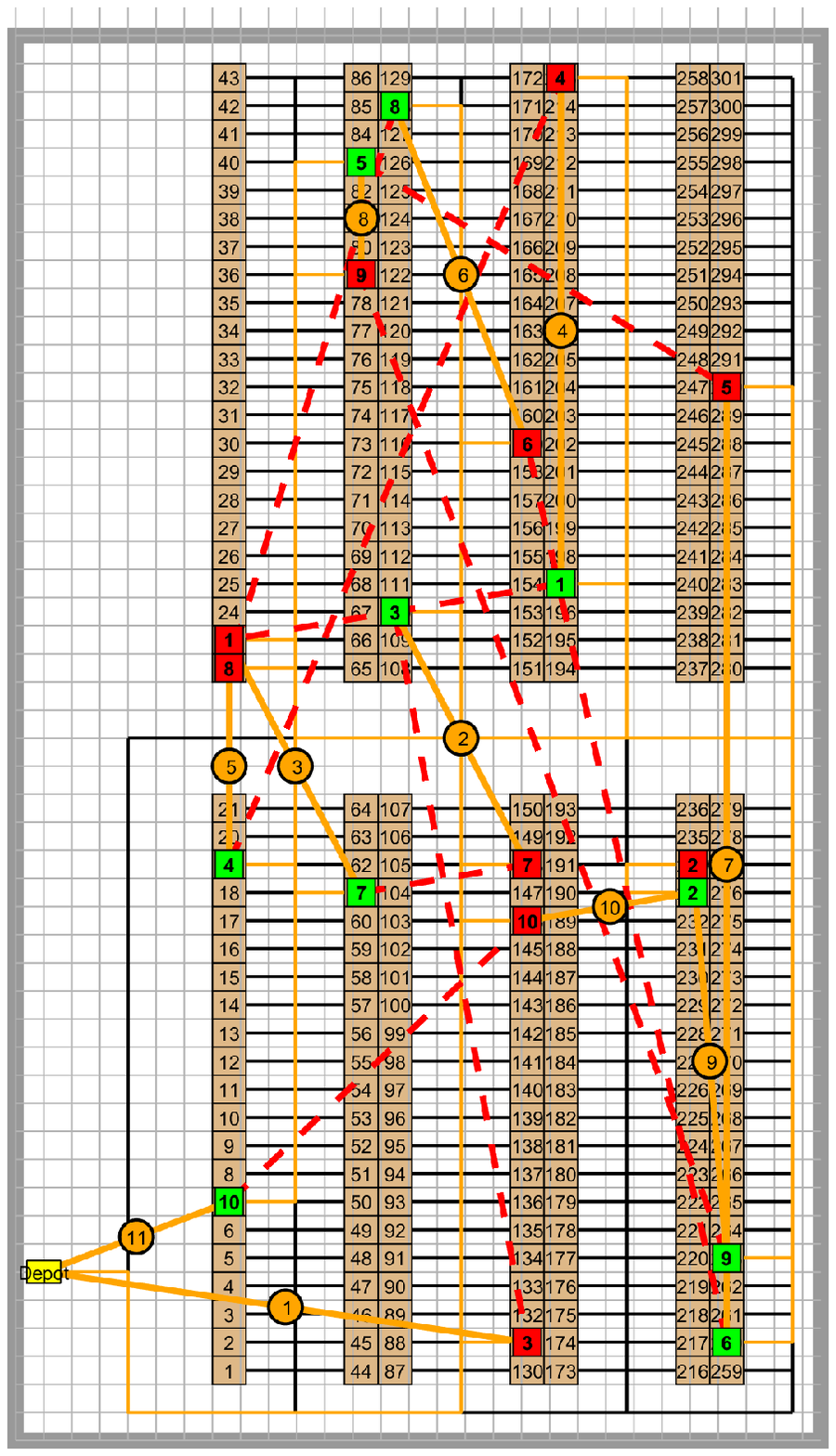} &
				\includegraphics[width=.5\textwidth, height=9cm, keepaspectratio]{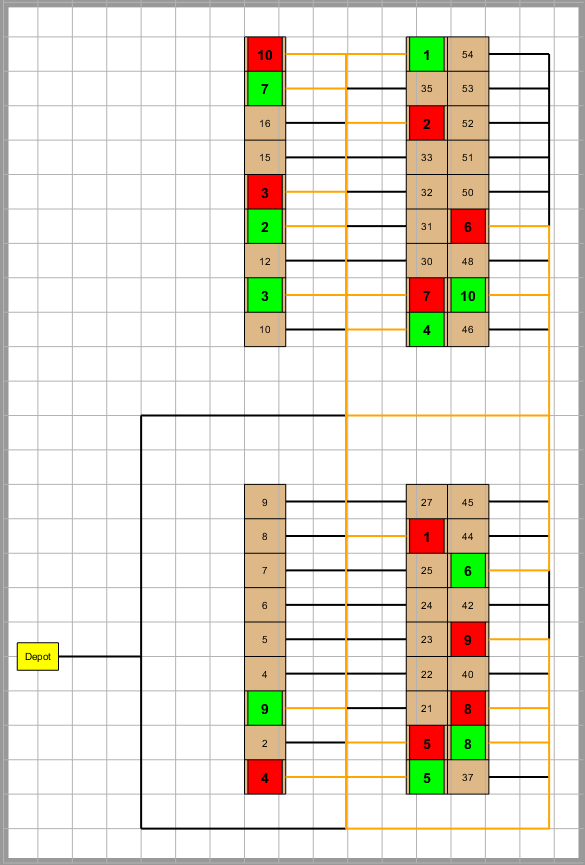}\\
				(a)&(b)\\
			\end{tabular}
			\caption{Warehouse (a) with static aTSP solution; (b) jobs for BD-mTSP.}%
			\label{fig:warehouse-depot}%
		\end{figure*}
		The truck can reach the central aisle and 
		has access to the south-side aisles from the depot. 
		The jobs are known in advance. 
		Hence, this requires a static asymmetric TSP to be solved. 
		The optimal solution derived with a binary programme is displayed within the figure. 
		This solution can be interpreted in the following way.
		The forklift leaves the depot and travels to storage location (SL) 131.
		Then the pallet is transferred to SL 110. 
		This completes transfer-job three.
		The forklift continues to job seven (SL 148 to 61).
		The entire tour in terms of transfer-jobs is  (0, 3, 7, 1, 4, 8, 6, 5, 9, 2, 10, 0), 
		where 0 represents the depot. 
	\end{example}

	\begin{example}[BD-mTSP AVH used in warehouse]
		Figure \ref{fig:warehouse-depot}b shows a warehouse network with $s=54$ storage locations,
		ten transfer-jobs $n=10$ .
		The shortest transfer job paths are highlighted in yellow.
		The distances of these are (17.5,  6.5,  6.5,  16.5,  4.5,  13.5,  9.5,  4.5,  18.4,  23.4).
		That means a minimum of 120.8 meters has to be travelled.
		Two trucks start and finish their tours at the depot (node zero).
		Balancing limits each truck to five jobs.
		Two jobs are done simultaneously by the two trucks. 
		At each discrete time step four jobs are revealed.
		That means, initially jobs 1 to 4 are known.
		From these two are executed. 
		The AVH selected jobs 1 and 4. 
		After the trucks completed these two jobs.
		Jobs 5 and 6 are available.
		This requires assigning the two trucks to two jobs from the set $\set{2,3,5,6}$.
		Job 2 and 3 are chosen. Next, jobs (5, 6, 7, 8) are visible.
		Job 5 and 7 are executed. The completion of these reveals the two final jobs.
		After doing them, the trucks return to the depots.
		In summary, tour 1 has the job sequence (0, 1, 2, 5, 8, 9, 0); and
		tour 2 (0, 4, 3, 7, 10, 6, 0).
		The distances travelled on these tours are 85.3 and 66.3 meters, respectively.
		Hence, the total travelled distance is 151.6 plus 120.8 meters. 
	\end{example}	

	The real-world instances are conceptually the same as the previous examples.
	However, the number of shelves per storage location increases,
	the size of the warehouse is larger, and the number of trucks increases.
	Zones and additional aisles may have to be considered.
	
	Table \ref{tab:whInstances} shows a summary of four warehouse instances.
	\input{whInstances.tex}
	The number of storage locations (\#SL) varies between 3,200 and 12,500.
	The number of jobs (\#jobs) per day ranges from 1,120 to 4,200 
	leading to different transfer-utilisations $\rho = \frac{\text{\#SL}}{\text{\#jobs}}$.
	The more storage locations the more forklift trucks $m$ are used.
	The vehicle-dynamics $D_a^m$ are between two and three, 
	which is reflected in the absolute dynamics $D_a$.
	The objective was to minimise the travelled distance $L$, 
	whilst balancing the number of jobs between forklift trucks.
	The distance between jobs is $L_j$ and the transfer-job distance is $L_i$.
	It is interesting to observe that the distances between jobs
	are a fraction (between 26.6\% and 39.0\%) of the transfers.
	The last column of the table states the total time the computations took in seconds.
	This varied between 507 and 1736 milliseconds, 
	demonstrating the suitability of the algorithm for online computations.

	\subsection{Taxi services \label{ssec:taxi}}
	
\begin{figure*}[htbp]\center
		\begin{tabular}{cc}
			\includegraphics[width=.4\textwidth, height=4cm, keepaspectratio]{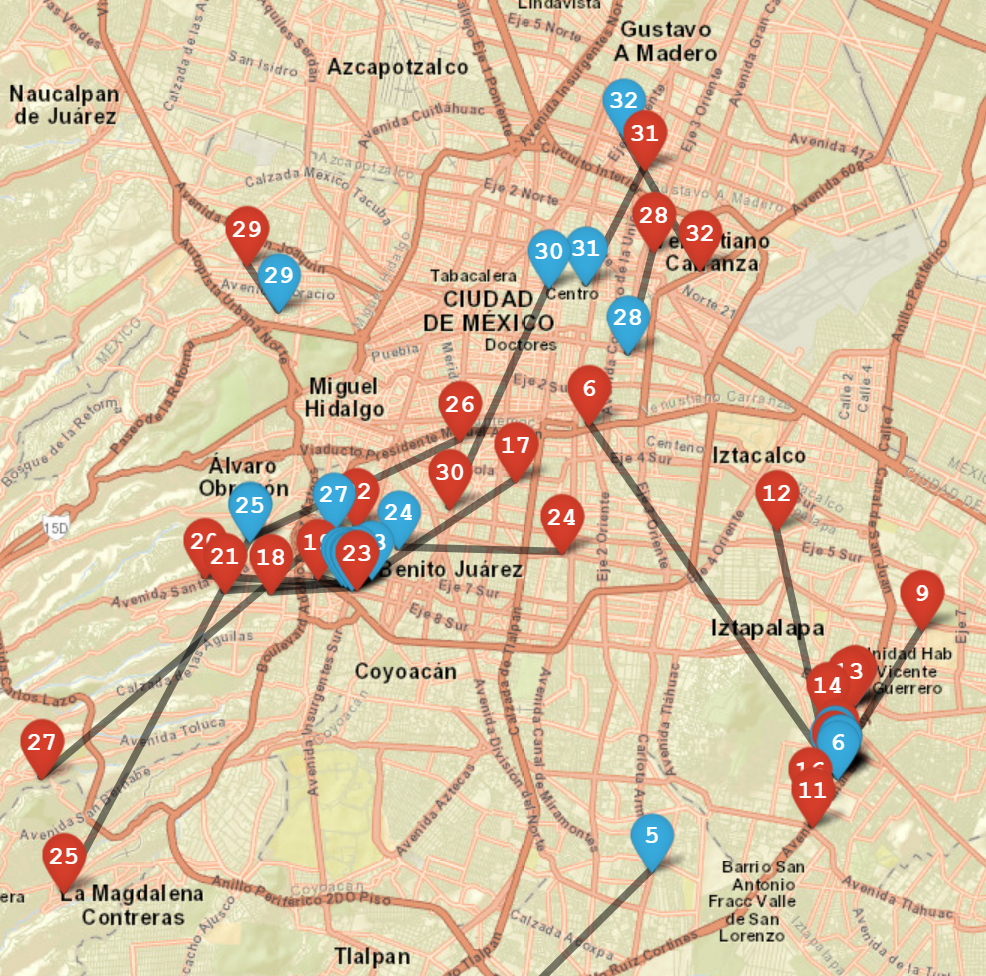}&
			\includegraphics[width=.6\textwidth, height=4cm, keepaspectratio]{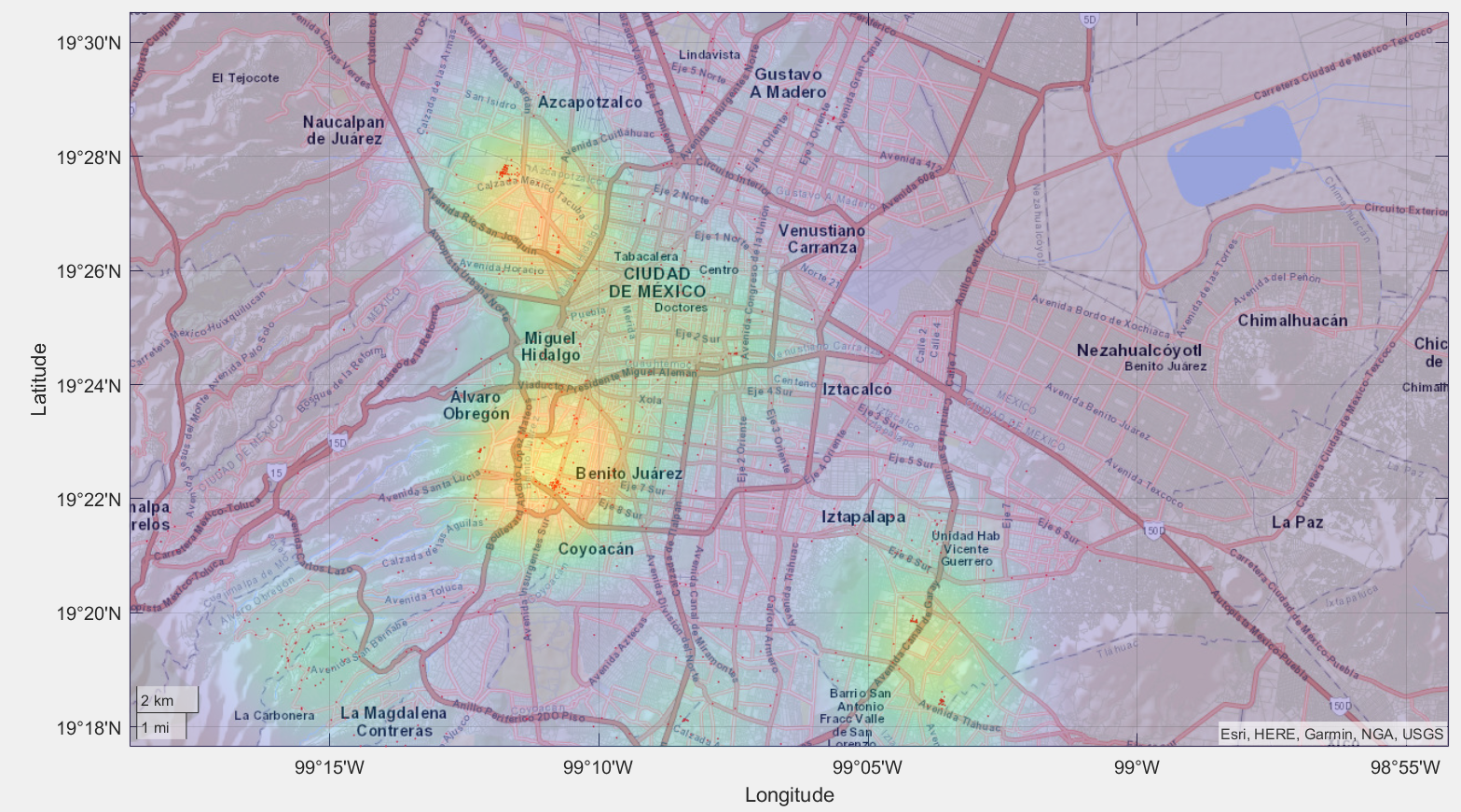}%
			\\ (a) & (b)
		\end{tabular}
		\caption{Mexico taxis (a) pick-ups and drop-offs on a day; (b) heatmap of pick-ups.}%
		\label{fig:taxi-mexico}%
	\end{figure*}

	Taxi data for Mexico city was obtained from Kaggle \cite{Navas2021}. 
	The data contains 12,694 records dating between 24/6/2016 and 20/7/2017.
	Importantly the records contain pick-up and drop-off coordinates and the corresponding times.
	On average, 33 daily journeys were recorded.
	Figure \ref{fig:taxi-mexico} (a) shows pick-ups (red) and drop-offs (blue) on 10th December 2016.
	The data was filtered by defining upper limits for 
	the waiting time (90 minutes), 
	trip duration (three hours), and 
	travelled distance (100 km).
	Trips were limited to the surrounding area of Mexico city.
	This was achieved by setting the pick-up latitudes to be between 19 and 20 degrees;
	and pick-up longitudes to be less than -98 degrees. 
	As a consequence, the dataset reduced to 9,571 records,
	which decreased the average of daily trips to 25. 
	The observed average distance is 6.65km.
	The average total distance of trips travelled per day is 166.6 km. 
	The travel duration is 27.0 min and waiting time is 8.72 min.
	Interestingly, the recorded time stamps are only between zero and the twelfth hour.
	
\begin{figure*}[htbp]\center
		\begin{tabular}{cc}
			\includegraphics[width=.5\textwidth, height=4cm, keepaspectratio]{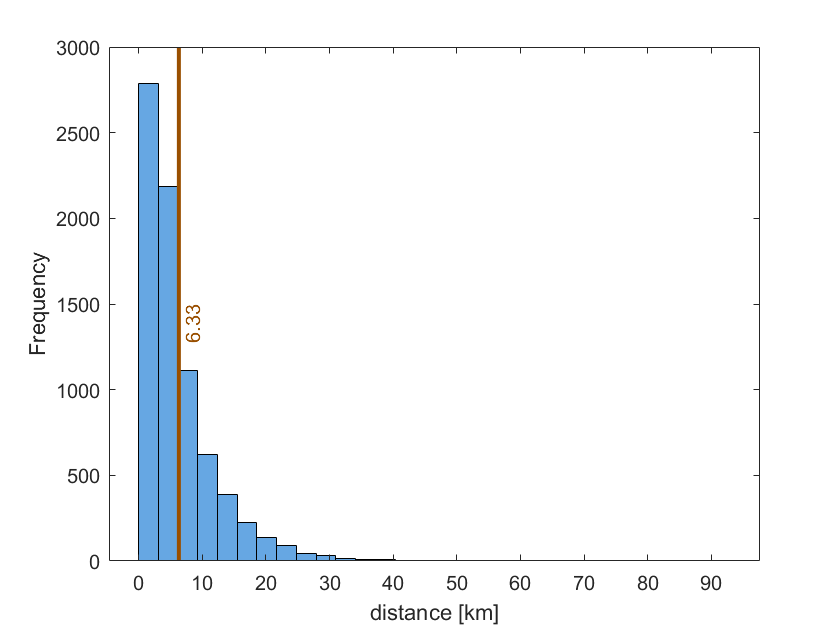}&
			\includegraphics[width=.5\textwidth, height=4cm, keepaspectratio]{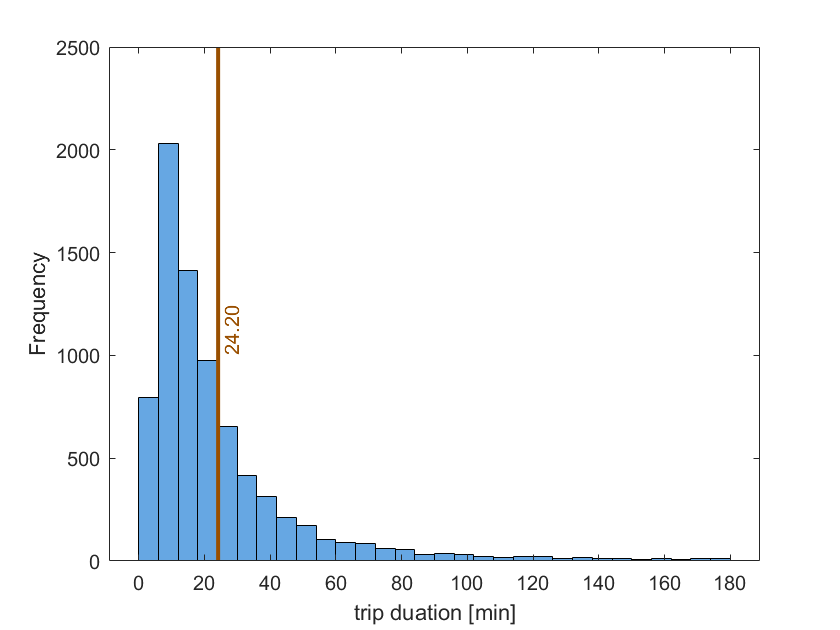}%
			\\ (a) & (b)
		\end{tabular}
		\caption{Mexico taxi histograms for (a) distances, and (b) trip duration.}%
		\label{fig:taxi-mexico-histo}%
	\end{figure*}

	Rather than finding the shortest route using the road network 
	the distance between a drop-off and a pick-up will be estimated using a detour factor.
	The detour factor will be derived from the given travelled distance data
	and using a geometric distance.
	The geometric distance between location 1 and 2 is determined using the Haversine formula.
	Locations are given by latitude $\lambda$ and longitude $\varphi$, which are in degrees.
	The degrees are converted to radians using $\pi = 180^{\circ}$.
	The latitude and longitude differences in radians are $d_\lambda = \frac{\pi}{180} (\lambda_2-\lambda_1)$ and
	 $d_\varphi = \frac{\pi}{180} (\varphi_2-\varphi_1)$.
	This allows us to define $\alpha$:
	\lequ{alpha}{\alpha = \sin^2 (\frac{d_\lambda}{2}) + 
		\cos {(\frac{\pi}{180}\lambda_1)}  \cos {(\frac{\pi}{180}\lambda_2)} 
	    \sin^2 {(\frac{d_\varphi}{2} }).}
    This allows us to compute the distance 
	\lequ{haversine}{d_h = 2r \arctantwo (\sqrt{\alpha}, \sqrt{1-\alpha}),}
	where $r$ is the earth radius estimated as 6,378.4 km.
	Equation \ref{eq:haversine} deviates from the simple formula 
	\lequ{simpledist}{d_s = \frac{\pi}{180} r \norm{\mat{\lambda_1-\lambda_2 \\ \varphi_1-\varphi_2}}_2}
	by about 80 meters on average (for distances in this dataset).
	For each trip, the recorded distance $d_j$ is given.
	This allows us to determine the detour factor for each trip.
	Trip distances where $\frac{d_j}{d_h}>3$ are specified as ``outliers''. 	
	Interestingly, these outliers constitute almost 20\% of the already filtered data.
	The remaining data was used to determine 
	the average detour factor 1.5752 for taxi trips in Mexico City
	(if the taxicab norm is used in Equation \ref{eq:simpledist} then the detour factor is 1.2148).
	The outlier trip-distances are repaired by using this detour factor.
	This gives an updated average travel distance for a trip of 6.33 km. 
	The average travel time is 24.2 minutes, and waiting time is 8.72 minutes.
	Figure \ref{fig:taxi-mexico-histo} show the corresponding distributions.
	25 trips divided by 12 hours gives 2.1 trip/hour.
	On average, a taxi can do two trips in 48.4 minutes not considering the time between drop-off and pick-up.
	Hence, as an initial estimate one taxi might be able to cover the entire workload.
	However, a closer inspection of the data shows that at least two taxis $m=2$ and sometimes more are required.
	Figure \ref{fig:taxi-mexico} (b) shows a heatmap of the pick-ups. 
	This allows the definition of a depot location $\lambda = 19.3702$ and $\varphi = -99.1799$.
	The distance matrix is derived between the drop-off and pick-up locations using Equation \ref{eq:haversine}.
	The first and last entry is the depot.
	
	The aim is to have a balanced workload between the taxis.
	Furthermore, the vehicle-dynamics are assumed to be $D_a^m \in \set{\frac{1}{2},1,2}$. 
	Representing the situations of (a) taxis waiting for customers needing a ride; 
	(b) a balanced situation; and (c) customers waiting for taxis to become available.
	
	The BD-mTSP AVH is used to determine the route-distances. 
	Before discussing the route and trip distances for all days an example of a single day is given.
	
	\begin{example}[Trips on a single day] \label{ex:taxi-single}
		On 10th of December 2016 there were 32 trips recorded in the Mexico taxi dataset (see Figure \ref{fig:taxi-mexico}a).
		After having applied the previously mentioned filters 29 trips remained.
		Five trip distance outliers ($ \frac{d_j}{d_h}>3 $) were observed (i.e. 17.2\%)
		and fixed using the detour factor 1.5752.
		This gives a total ``internal'' (from pick-up to drop-off) trip-distance of 152.5km.
		Two taxis ($m=2$) operate to fulfil the trips.
		Table \ref{tab:taxiDec10} shows the total distances travelled 
		depending on vehicle dynamics $\frac{1}{2}$, 1 and 2.
		\input{taxiDec10} 
		The distances between drop-offs and pick-ups are determined using the BD-mTSP AVH.
		We observe that the distances to reach customers exceeds the ones for the passenger trips.
	\end{example}
	
	The Mexico data was analysed for the entire period (392 days).
	However, only days with more than ten trips per day were used.
	This constraint reduced the dataset to 327 days.
	The average travelled daily trip-distance is $\overline{L_i} = 149.3$ km 
	with a standard deviation of $\sigma_{L_i} =77.2$ km.
	
	Table \ref{tab:taxi-analysis} shows results from analysing the Mexico City data.
	\input{tab_taxi-analysis}
	It can be seen (last column) that the average distance travelled without passengers ($L_j$)
	exceeds the average daily trip-distance ($L_i$).
	It is also observed that the daily travelled distances ($L_i$, $L_j$)  vary substantially.
	
	Considering the vehicle-dynamics $D_a^m$ 
	we can confirm the validity of the distance-vehicle-dynamic proposition empirically.
	That means, the route-distances $L_j$ are highest at $D_a^m = 1$ and are lower for $D_a^m \in \set{\frac{1}{2},2}$.


%% file: whInstances.tex
\begin{table*}[htbp]
	\centering
	\caption{Warehouse instances}
\begin{tabular}{ccccc|cccc|c}
	\toprule
	\#SL & \#jobs & $\rho$ & $m$  & $D_a$ & $L$  & $L_j$ & $L_i$ & $L_j/L_i$ & sec \\
	\midrule
	3,200  &        1,120  & 35.0\% &             18  &             36  &   142,783  &           40,059  &        102,724  & 39.0\% &        0.507  \\
	5,600  &        2,700  & 48.2\% &             27  &             81  &   413,322  &           94,683  &        318,639  & 29.7\% &        0.894  \\
	7,300  &        3,500  & 47.9\% &             30  &             60  &   638,265  &        162,326  &        475,939  & 34.1\% &        0.962  \\
	12,500  &        4,200  & 33.6\% &             44  &           132  &   933,045  &        196,101  &        736,944  & 26.6\% &        1.736  \\
	\bottomrule
\end{tabular}%

	\label{tab:whInstances}%
\end{table*}%

%% file: taxiDec10.tex
\begin{table}[htbp]
	\centering
	\caption{Total tour lengths for December 10th.}
	\begin{tabular}{r|rrrr}
		\toprule
		\multicolumn{1}{c|}{$Da$} & \multicolumn{1}{c}{$L$} & \multicolumn{1}{c}{$Lj$} & \multicolumn{1}{c}{$Li$} & \multicolumn{1}{c}{$Lj/Li$} \\
		\midrule
		1    & 368.5 & 216  & 152.5 & 141.6\% \\
		2    & 370.7 & 218.2 & 152.5 & 143.1\% \\
		4    & 332.5 & 180  & 152.5 & 118.0\% \\
		\bottomrule
	\end{tabular}%
	\label{tab:taxiDec10}%
\end{table}%

%% file: tab_taxi-analysis.tex
\begin{table}[htbp]
	\centering
	\caption{Analysis results of taxis in Mexico.}

\adjustbox{width=\columnwidth}{ 
\begin{tabular}{ccc|ccc}
	\toprule
	$m$  & $D_a^m$ & $D_a$ & \ms{L} & \ms{Lj} & \ms{L_j/L_i} \\
	\midrule
	\multirow{3}[2]{*}{2} & $\frac{1}{2}$ & 1    & 383.0±214.5 & 233.7±146.1 & 160.3\%±54.9\% \\
	& 1    & 2    & 414.1±237.9 & 264.8±170.4 & 180.9\%±66.1\% \\
	& 2    & 4    & 349.1±184.8 & 199.8±116.4 & 139.6\%±49.7\% \\
	\midrule
	\multirow{3}[2]{*}{3} & $\frac{1}{2}$ & $\frac{3}{2}$ & 361.8±190.8 & 212.6±122.5 & 148.9\%±51.7\% \\
	& 1    & 3    & 404.8±217.2 & 255.5±149.7 & 178.7\%±66.7\% \\
	& 2    & 6    & 342.9±168.1 & 193.6±101.0 & 139.9\%±55.7\% \\
	\midrule
	\multirow{3}[2]{*}{4} & $\frac{1}{2}$ & 2    & 365.0±176.5 & 215.7±108.2 & 155.4\%±57.9\% \\
	& 1    & 4    & 400.6±203.7 & 251.3±137.1 & 178.8\%±68.4\% \\
	& 2    & 8    & 348.6±160.6 & 199.3±94.8 & 146.8\%±61.8\% \\
	\bottomrule
\end{tabular}%
}

	\label{tab:taxi-analysis}%
\end{table}%

%% file: sec-CAM.tex
\section{Continuous Approximation Model}\label{sec:CAM}

In this section we will give the functional relationship
that relates the solution distance 
to the number of salesmen (vehicles), customers (nodes) and sequential dynamics.

\subsection{Related Work}
One of the first researches addressing the above for the static TSP 
was the work by \cite{beardwood1959shortest} called a shortest path through many points.
This led to the famous Beardwood-Halton-Hammersley theorem,
which gives an asymptotic formula for the length of a TSP route.
\begin{theorem}[The Beardwood–Halton–Hammersley (BHH) Theorem]
	Let $\left\{X_{1}, \ldots, X_{n}\right\}, n \geq 1,$ be a set of random variables
	in $\mathbb{R}^{d},$ independently and identically distributed with bounded support. 
	Then the length $L_{n}$ of a shortest TSP tour through the points $X_{i}$ satisfies
	\lequ{BHH-theorem}{\begin{split}
		L_{n} / n^{(d-1) / d} \rightarrow \beta_{d} \int_{\mathbb{R}^{d}} f(x)^{(d-1) / d} d x, \\
		\text {with probability 1, as } n \rightarrow \infty \end{split}}
	where $f(x)$ is the absolutely continuous part of the distribution of $X_{i}$ 
	and $\beta_{d}$	depends on $d$ but not on its distribution.
\end{theorem}
This theorem influenced probability theory, physics, computer science and operational research.
The case when $d=2$ is of interest to us, i.e.
a TSP length is `almost always' asymptotically proportional to $\sqrt{nv}$,
where $n$ represent the number of points in a bounded plane region with area $v$.
The exact value of $\beta_2$ is still unknown \citep[p23]{applegate2006traveling}.
It was approximated as 0.7313 for $n=1,000$.
\citet[p497ff]{applegate2006traveling} provide more estimates for $\beta_2$
from various authors. 
The $\beta_2$ values are between 0.7765 and 0.7241 for nodes between $n=100$ and $n=2,500$.
Most approaches to find a function $L_d = f(X)$ for the TSP use heuristics or regression models. 
Building on the above work
\cite{ccavdar2015distribution} gave a distribution-free TSP tour length estimation model for random graphs.
Their approach used a regression model based on sampling probability distributions.

An interesting empirical formula for the CVRP was developed by
\cite{Eilon1971} $(L / N) \cong 1.8 \bar{\rho}[(1 / C)+(1 / \sqrt{N})]$,
where the depot is in the centre and $N$ points are uniformly distributed in a square.
Here, $1.8 \bar{\rho}/C$ represents the general location to reach a point and
$1.8 \bar{\rho}/\sqrt{N}$ the detour distance. 
$\bar{\rho}\cong 0.382 A^{1 / 2}$, where $A$ is the area of interest,
and $C$ is the number of items. 
This was the starting ``stone'' for \cite{daganzo1984distance} ,
who consider a special case of the CVRP, which agrees with the mTSP definitions given earlier.
They derived a formula for the CVRP by extending the above TSP findings.
This was achieved by cluster-first and route-second, i.e. a heuristic model.
They emphasise the importance of choice of shape (e.g. slenderness) for the clusters.
Their length formula was expressed as: 
\lequ{daganzo}{L(\mathscr{A}) \cong 1.27 \times\left(2(\bar{\rho} / C)+0.57 \delta^{-1 / 2}\right),}
where $\bar{\rho}$ is the average Euclidean distance of the node locations to the depot,
$C$ the number of points in an sector, 
$\mathscr{A}$ the total area, and
$\delta$ is the density of the area.
In \cite[p5,eq14]{garn2020closed} I proposed a continuous approximation model for the balanced static mTSP:
\lequ{final-mTSP}{L(\mathscr{A}) \approx   138.2 n^{0.44}  + 88.1(m-2),}
where $\mathscr{A}$ is a discrete 100 by 100 area and customers are uniformly distributed on the grid.

\cite{franceschetti2017continuous} review literature on continuous approximation models in freight distribution management.
They identified \cite{beardwood1959shortest} and \cite{daganzo1984distance} as seminal works (see above for details).
Their concise review contains several more formulations and applications of continuous approximation models 
containing several valuable ideas and formulations.
However, their work did not identify an approximation formula for the dynamic-mTSP, 
which indicates a potential gap in the body-of-knowledge.
Similarly, \cite{ansari2018advancements} discuss advancements in continuous approximation models for logistics and transportation systems by focusing their review on literature between 1996 and 2016. 
Their work covers a fast range of applications and problem instances.
This confirms that approximation models for the dynamic mTSP are at best sparse.
Nevertheless, there is one report by \cite{erera2003dynamic},
which introduces a ``threshold global sharing'' scheme utilising a real-time re-optimisation control for the VRP.
At first it appears to have similarities with the CVH.
However, it differs because of an initial partitioning approach, which is typical for Daganzo's studies.
Roughly speaking, their model is based on a spherical area. 
Customers' density, expected demand, and standard deviation are used.
Vehicles' capacity is set. Additionally, a buffer factor $\alpha$ is given 
which defines the number of standard deviations for the total customer demand.
Overall, their approximation model has a high level of complexity, 
which makes it difficult to use for verification purposes.

In the next section I will propose a generic approach to derive CAM models.

\input{SeqDyn}

%% file: SeqDyn.tex
\subsection{Continuous BD-mTSP approximation}
In previous sections we introduced 
and used two algorithms to determine routes and associated lengths.
This section will show how to create a formula that can estimate route lengths.

\begin{definition}[real-world model]
	A \textit{real-world model} $f$ (systematic information) is explained by:
	\equ{y=f(X)+\epsilon,}
	where $y$ is the output (response),
	$X = (x_1, \cdots, x_p)$ is the input (known also as predictors, variables, features)
	and $\epsilon$ is the random error.
	All inputs and outputs are real with $X \in \real^{n\times p}$ and $y \in \real^{n}$.
	Here, $n$ is the number of observations and $p$ is the number of features.
\end{definition}

As a first interpretation, $y$ represents the route length obtained
from the BD-AVH algorithm's route output.
The algorithm itself represents the function $f$.
The algorithm's original inputs are a distance matrix $D$, vehicles $m$ and absolute dynamics $d$.
The distance matrix can be obtained from the number of customers $n-1$ as will be explained below.

As a second interpretation, which will be used later, $y$ represents the average route length.
This average is obtained from running the BD-AVH algorithm multiple times. 
Consequently, the BD-AVH will be called \textit{average BD-AVH}.
The output values of $y$ are also known as the \textit{observed values}.
The second interpretation allows to get statistical stability.

\begin{definition}[prediction model]
	A prediction model $\hat{f}$ (estimate for $f$) generates predictions $\hat{y}$ 
	using existing (real) input $X=\mat{X_i}=\mat{x_j}$ (rows $X_i$, columns $x_j$):
	\equ{\hat{y}=\hat{f}(X).}
	$\hat{f}$ is influenced by the degree of \textit{flexibility}, 
	e.g. used features.
\end{definition}

Here, $ \hat{f} $ is the continuous approximation model (CAM) for the average BD-AVH.
The purpose of the remainder of the section is to find this CAM.
The predictions will be called \textit{modelled values}.
The prediction model's objective is to estimate the route length $\hat{y}$
given the number of vehicles $m$, the number of customers $n-1$ 
and the absolute dynamics $D_a$ such that the root mean squared error is minimised.

The quality of the fitted $\hat{f}$ can be described using the root mean squared error (RMSE):
\equ{\text{RMSE}=\frac{1}{n} \sqrt{\sum_{i=1}^{n}\left(y_{i}-\hat{f}\left(X_{i}\right)\right)^{2}},}
where $y\in \real^n$ and $X_i \in \real^p$ is the $i$-th row.


The idea of the proposed approach is to run the BD-AVH on uniformly distributed customers in the Euclidean plane.
A range of fleet sizes and sequential absolute vehicle-dynamics are provided.
Regressing on the route lengths and feature matrices will provide the approximation formula.

The analysis will be restricted to the unit square area.
The advantage of this is the scalability to other area sizes.
Within this area $n$ nodes with an $x$ and $y$ coordinate are created.
These are uniformly distributed numbers within the interval $ (0,1) $.
The first node represents the depot. 
That means, all customer locations including the depot are randomly located within $\mathscr{A} = (0,1)^2$.
The above explains, why $n$ can be used as input for the BD-AVH instead of the distance matrix $D$.

There are alternatives to the placement of customers and depots, which future work can consider.
For instance, \cite{uchoa2017new} generated CVRP benchmark instances 
by placing the depot randomly (same as this study's assumption), central (centre of the grid) or eccentric (corner). 
Their work suggests random (same as this work), clustered and random-clustered positioning of customers.
Hence, there are nine combinations (scenarios) of depots and customer locations that could be considered.
Depending on the application of the CAM other distribution estimates may be considered
such a Normal distribution (with the depot being in the mean-centre), Weibull distribution, Logistic distribution or Gamma distribution.
This might better reflect some real-world scenarios. 
For instance, the drop-off locations for taxi services (Section \ref{ssec:taxi}) are better approximated as normal distribution 
than a uniform distribution.
When modelling warehouse locations (Section \ref{ssec:wh}),
uniformly distributed customer positions on integer grids are more appropriate.
Generally, as an initial step, the distribution of the customers and placement of the depot need to be determined.
For this purpose, distribution fitting methods \cite[p175]{Garn2018} can be used.
The framework introduced below can be applied to obtain a set of CAMs for the above mentioned scenarios.
However, this work focuses on placing the depot and customers using a uniform random distribution.

\begin{figure*}[htbp]
	\centering
	\begin{tabular}{ccc}
		\includegraphics[width=.3\textwidth, height=6cm, keepaspectratio]{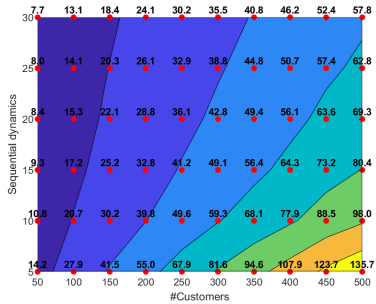}&
		\includegraphics[width=.3\textwidth, height=6cm, keepaspectratio]{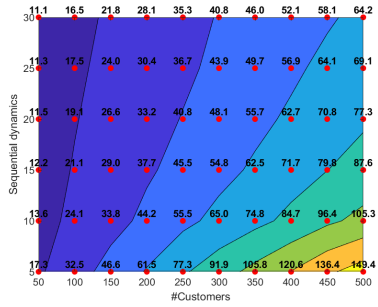}&
		\includegraphics[width=.3\textwidth, height=6cm, keepaspectratio]{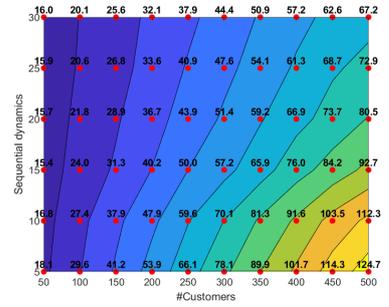}\\
		(a) $m=1$ & (b) $m=4$ & (c) $m=7$
	\end{tabular}
	\caption{Average distances given sequential absolute dynamics in relation to nodes and vehicles for the BD-AVH. }
	\label{fig:DAVH-uniform}
\end{figure*}

Figure \ref{fig:DAVH-uniform} shows 
the average distances found by the BD-AVH using absolute sequential dynamics.
The results seen in the figure motivated to choose the BD-AVH rather than the BD-CVH, 
because route lengths are shorter by 1.61\%. 
This
 percentage was obtained by determining the average relative difference between BD-AVH and BD-CVH
for all uniform test instance solution values.


It is interesting to note the upper bound $(n+1) \sqrt{2}$, if sequential dynamic visibility is one 
and customers are always in opposite corners.
The instance $m=1, d=5, n=500$ has the upper bound 708.5 and BD-AVH distance 135.7 (see Figure \ref{fig:DAVH-uniform}).
The lower bound is trivial.

The average BD-AVH algorithm uses 30 test instances for each configuration.
A configuration $c=\mat{m,n,d}$ is a 3-tuple with:
$$ m\in \set{1,\dots,7}, n \in \set{50, 100,\dots,500}, d\in \set{5,10,\dots,30}.$$
Hence, the total number of configurations is $7 \times 10 \times 6 = 420$. 
The configurations are created by looping through the three sets in the sequence given above.
They are captured in $X\in \setN^{420\times 3}$.
\begin{example}[Configuration 129]
	The first configuration is $X_{1}=\mat{m=1,n=50,d=5}$.
	Iterating through the absolute dynamics set we get
	$X_{2}=\mat{m=1,n=50,d=10}$ and $X_{3}=\mat{m=1,n=50,d=15}$.
	Once six iterations are done the next element in $n$ is used, 
	i.e. $X_{7}=\mat{m=1,n=100,d=5}$.
	That means, there are 60 iterations before $m$ is incremented.
	$X_{129}=\mat{m=3,n=100,d=15}$ is the 129$^\text{th}$ configuration having 3 vehicles, 100 customers and sequential dynamics of 15.
	Here, 100 customer locations are randomly created 30 times. 
	Leading to 30 distance values, where $y_{129}=19.8$ is the average distance.
\end{example}

We will derive an approximation formula (model) using multivariate regression, i.e.
\lequ{LR}{y = X b+\epsilon \ie y \perp X \ie b = (\tp{X}X)^{-1}  \tp{X}  y.}
Here, $y$ is the total travelled average distance and $X$ are features.
$y$ represents ``observed'' values using the BD-AVH algorithm.
In our case the BD-AVH algorithm was run 30 times and its average saved in $y_i$.
Since there are 420 configurations $y \in \real^{420}$.
The number of features will vary depending on our approach, 
but each feature $x_j$ will have 420 elements.
The prediction model is given by:
\lequ{mLR}{\hat{y} = \hat{f}(X) = X b.}

\begin{figure*}[htbp]
	\centering
	\begin{tabular}{cc}
		\includegraphics[width=.45\textwidth, height=6cm, keepaspectratio]{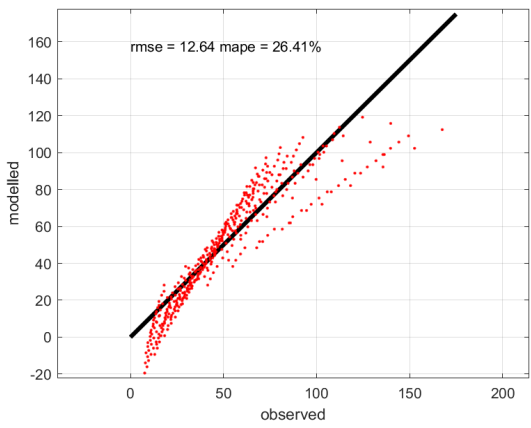}&
		\includegraphics[width=.45\textwidth, height=6cm, keepaspectratio]{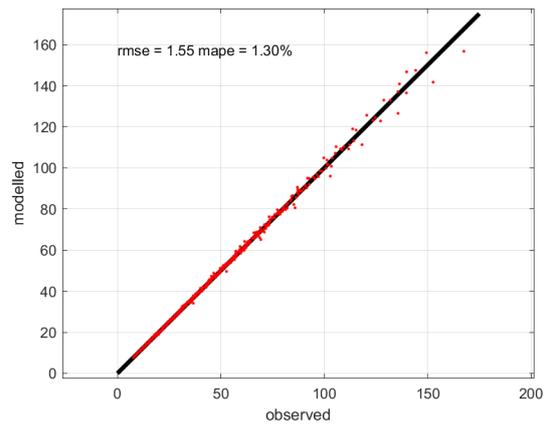}\\
		(a)  & (b)
	\end{tabular}
	\caption{Distances modelled with linear regression compared to observed ones from the BD-AVH 
		using (a) features as-is, (b) feature map. }
	\label{fig:CAM-LR}
\end{figure*}

A naive approach is to use the features as-is.
That means the number of vehicles $x_1$, customers $x_2$ and sequential absolute dynamics $x_3$.
This defines $X=\mat{x_1&x_2&x_3}$. 
I will refer to these as the \textit{base-features} (or \textit{feature basis}),
i.e. $X$ consists out of three base-feature vectors.
That means, $X$ represents the 420 configurations discussed earlier.
Using Equation \ref{eq:LR} allows us to derive the model.
This in turn allows us to compute all modelled values.
Figure \ref{fig:CAM-LR} (a) compares the observed and modelled distances with each other.
In a perfect model the points (observed, modelled) lay on the black line, i.e.
the modelled distances are the same as the observed ones.
However, it can be seen this is not the case.
To get a better idea about the accuracy of the prediction method another quality measure is introduced.
The \textit{mean absolute percentage error} (MAPE) focuses on relative deviations and is commonly found in forecasting:
\lequ{mape}{\text{MAPE} = \frac{1}{n} \sum_{i=1}^n\left|\frac{y_i-\hat{y}_i}{y_i}\right|,}
where $\hat{y}_i = \hat{f}\left(X_{i}\right)$.
The main issue of the MAPE is its sensitivity with small values, and that it is not defined when any $y_i=0$.
However, the advantage is its explainability.
For the initial model, the RMSE is 12.64 and the MAPE is 26.41\% which are improvable.

We know from our previous and related work 
that the distance is proportional to the square root when dynamics are disregarded.
Additionally, from Figure \ref{fig:DAVH-uniform} there is polynomial behaviour visible.
Hence, new features are engineered from the base-features using combinations of powers.
The powers are $p \in \set{0,\frac{1}{2},1,2}$ 
\footnote{Powers $p \in \set{0,1,2,3}$ perform equally well: RMSE 1.67, MAPE 1.69\%. }.
The resulting feature map consists out of 64 features:
$$ x_1^0 x_2^0 x_3^0,  ~x_1^0 x_2^0 \sqrt{x_3},~x_1^0\sqrt{x_2}x_3^0,\dots,~ x_1^1 x_2^2 x_3^2, ~x_1^2 x_2^2 x_3^2.$$
This increases the flexibility of the model.
By the way, the base-features are contained in the map.
Figure \ref{fig:CAM-LR} (b) compares the modelled and observed distances.
The RMSE is 1.55, which is low considering that $\bar{y}$ is 51.9.
The MAPE is 1.30\% which means that the model is very well fitted or overfitted.

As a technical side note, the computation of inverse of $(\tp{X}X)$ returns a matrix that is close to singularity or badly scaled 
(absolute value magnitudes range between $10^{-15}$ and $10^6$).
Using the pseudo-inverse instead leads to an absolute magnitude range between $10^{-10}$ and $1$.
However, the RMSE and MAPE increase to 4.55 and 7.41\% respectively.

\begin{figure*}[tbph]
	\centering
	\includegraphics[width=0.9\linewidth]{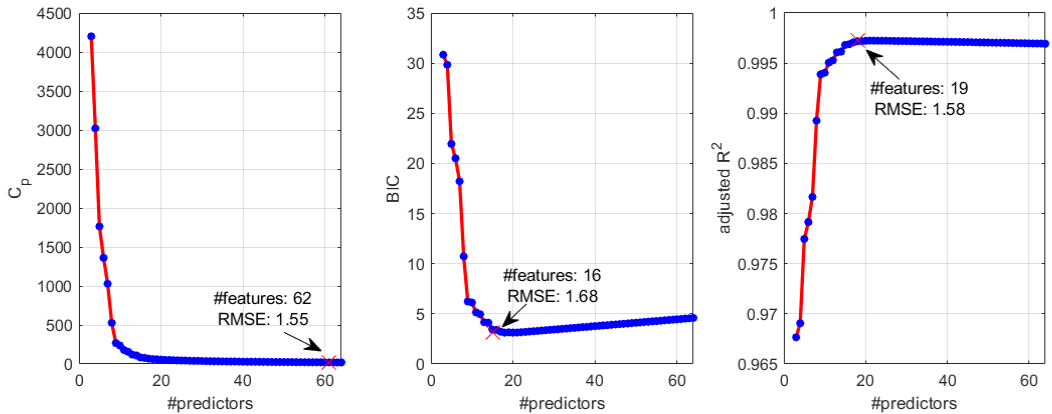}
	\caption{$C_p$, BIC and adjusted $R^2$ are shown for the best model for each feature step.}
	\label{fig:subsetselection}
\end{figure*}

In order to improve the interpretability of the model and avoid over-fitting - regularisation can be applied.
Three common approaches are subset selection, shrinkage, and dimension reduction \citep{james2013introduction}.
The best-subset selection approach requires the fitting of $2^{64}$ regression models.
This is computationally infeasible. 
There are two alternatives: forward or backward stepwise selection,
which are greedy search algorithms. 
Both of them are computationally efficient and have similar results.
The backward stepwise selection algorithm was used - see \cite{james2013introduction} for details.
This returns 64 models with features ranging between 1 and 64.
Choosing the optimal model can be achieved using  Colin Lingwood Mallows's cross-validation prediction criteria $C_p$, 
Akaike information criteria (AIC), 
Bayesian information criteria (BIC) or adjusted-$R^2$.
Figure \ref{fig:subsetselection} shows $C_p$, BIC and adjusted $R^2$ for the best model for each feature step.
Models with 0, 1 and 2 features were omitted in the display due to their large criteria values and large RMSEs ($>6.13$). 
The cross-validation prediction errors $C_p$ suggest a model with 62 features,
which appears to be a poor choice considering that the RMSE changes only slight 
between feature 20 (rmse: 1.56, mape: 1.39\%) and 62 (rmse: 1.55 , mape: 1.27\%).
The case with three selected features is of interest,
because this is the same number of features as in the naive approach but with a RMSE of 5.44 instead of 12.64 (MAPE: 9.83\% instead of 26.41\%).
Recall the three base-features are vehicles $x_1$, customers $x_2$ and dynamics $x_3$. 
The selected features (as polynomials of the base-features) for our first model are:\\
\begin{center}
	\begin{tabular}{cc} \hline
		polynomial & 	b \\ \hline
		$x_2$ & 	0.39051 \\
		$x_2 \sqrt{x_3}$ & 	-0.05477 \\
		$x_1 x_2 x_3$ & 	0.00023 \\
		\hline
	\end{tabular}
\end{center}
An advantage is the ease of using this model via:
\lequ{CAM-3f}{ L = 0.391x_2-0.055 x_2 \sqrt{x_3}+2.33\cdot 10^{-4} x_1 x_2 x_3.}
\begin{example}[Predicted distance with CAM]
	Let $\mat{x_1& x_2& x_3}=\mat{m=3,n=100,d=15}$ be an input configuration.
	The observed average distance is $y=19.8$, which was computed using the BD-AVH.
	The predicted average distance \ref{eq:CAM-3f} is $\hat{y}=L=18.85$ obtained from Equation (\ref{eq:CAM-3f}).
	The corresponding percentage error is 4.8\%.
\end{example} 
We observe that $x_2$ is in all terms of (\ref{eq:CAM-3f}). 
This emphasises the importance of the total number of nodes (customers). 
However, unlike to the static mTSP the length is not proportional to the square root of customers.
A reason could be that the visibility of the customers 
is limited by sequential dynamics.
However, the square root is applied to the sequential dynamics scope $x_3$ base-feature.
This makes sense, since this represents the visible nodes, 
which are used when assigning the vehicles to the customers.
The number of vehicles is found in the $x_1 x_2 x_3$ term.

The recommended model according to BIC derived from the backward selection algorithm is the model with 16 features and a RMSE of 1.68 (MAPE: 2.11\%) is shown.
Table \ref{tab:cam-f16} shows the polynomial-features and model coefficients.
\input{table-cam-f16}
This (second) model extends the previously mentioned model with 13 new features.
$\sqrt{x_3}$ can be found in five terms confirming the importance of sequential dynamics.
The total number of nodes (customers) $x_2$ is present in almost all terms with the exception of term 7.
The model suggested by the BIC is less accurate than the ones recommended by the other information criteria
but has less bias. 
In this case, it appears that the BIC is most effective in reducing the number of features.

However, an even better compromise between model quality and interpretability 
appears to be model three with 9 selected features (see Table \ref{tab:CAM-9f}).
\begin{table}[htbp]
	\caption{BD-mTSP CAM with 9 features.}	
	\label{tab:CAM-9f}
	\begin{center}
	\begin{tabular}{cc} \hline
		polynomial & 	b \\ \hline
		$\sqrt{x_1}  \sqrt{x_3}$ & 	0.52829 \\
		$x_2$ & 	0.29958 \\
		$x_1 x_2$ & 	0.17818 \\
		$x_1 x_2 \sqrt{x_3}$ & 	-0.08168 \\
		$x_2 \sqrt{x_3}$ & 	-0.03651 \\
		$x_1^2 x_2$ & 	-0.02354 \\
		$x_1^2 x_2 \sqrt{x_3}$ & 	0.01102 \\
		$x_1 x_2 x_3$ & 	0.00927 \\
		$x_1^2 x_2 x_3$ & 	-0.00126 \\
		\hline
	\end{tabular}
	\end{center}
\end{table}
This model has a RMSE of 2.35 and a MAPE of 2.82\%.
The model is acceptable by all information criteria considering Figure \ref{fig:subsetselection}.

The CAM models can be used to approximate the expected distance when nodes are uniformly distributed in the unit square area.
The three proposed models vary in flexibility and accuracy. 
The first model is simple to use (three features) and has a MAPE of 9.83\%.
The second model has an accuracy of 2.11\% (MAPE) but requires 16 features.
The third model is a compromise using 9 features and having an accuracy of 2.82\%.

%% file: table-cam-f16.tex
\begin{table}[htbp]
	\centering
	\caption{Continuous approximation model using 16 features.}
	\adjustbox{width=\columnwidth}{
	\begin{tabular}{ccc|ccc}
		\toprule
		\tiny{term} &           polynomial           & \small{    b    } & \tiny{term} &        polynomial         & \small{    b    } \\ \midrule
		\tiny{ 1  } &      $\sqrt{x_1} x_2$       & \small{- 2.82526} & \tiny{ 9  } &     $x_1 x_2 x_3$      & \small{ 0.11576 } \\
		\tiny{ 2  } &            $x_2$            & \small{ 1.93401 } & \tiny{ 10 } &      $x_1^2 x_2$       & \small{- 0.10610} \\
		\tiny{ 3  } &          $x_1 x_2$          & \small{ 1.56537 } & \tiny{ 11 } &       $x_2 x_3$        & \small{ 0.08903 } \\
		\tiny{ 4  } & $\sqrt{x_1} x_2 \sqrt{x_3}$ & \small{ 1.40787 } & \tiny{ 12 } & $x_1^2 x_2 \sqrt{x_3}$ & \small{ 0.06008 } \\
		\tiny{ 5  } &    $x_1 x_2 \sqrt{x_3}$     & \small{- 0.82816} & \tiny{ 13 } &    $x_1^2 x_2 x_3$     & \small{- 0.00922} \\
		\tiny{ 6  } &      $x_2 \sqrt{x_3}$       & \small{- 0.81389} & \tiny{ 14 } &    $x_1 x_2 x_3^2$     & \small{- 0.00053} \\
		\tiny{ 7  } &  $\sqrt{x_1}  \sqrt{x_3}$   & \small{ 0.52925 } & \tiny{ 15 } & $\sqrt{x_1} x_2 x_3^2$ & \small{ 0.00041 } \\
		\tiny{ 8  } &    $\sqrt{x_1} x_2 x_3$     & \small{- 0.17718} & \tiny{ 16 } &   $x_1^2 x_2 x_3^2$    & \small{ 0.00006 } \\ \bottomrule
	\end{tabular}%
	}
	\label{tab:cam-f16}%
\end{table}%

%% file: sec-Conclusion.tex
\section{Conclusion}\label{sec:Conclusion}
Dynamic routing is essential in many real-world scenarios.
This work focused on dynamics for the balanced mTSP.
However, this work can be easily applied to the capacitated vehicle routing problem (CVRP).
It is closely related to dial-a-ride problem (DARP) as explained in Section \ref{ssec:dynamics-scope}.
Several types and scopes of dynamics were proposed.
This work focused on the sequential dynamic scope.
Future work can investigate variable dynamics.
Two algorithms: balanced-dynamic closest vehicle heuristic (BD-CVH) 
and balanced-dynamic assignment vehicle heuristic (BD-AVH)
were developed.
It would be interesting to compare the BD-AVH to insertion-heuristics.
Several test instances (including derivates from the TSPLIB) 
gave insights about the behaviour of dynamic routing.
The proposed test instances and solutions can be used
as benchmark reference.
The observed distances indicated that in general the BD-AVH 
algorithm leads to slightly better solutions (about 3\%).
A comparison to static routing on these test instances 
gave an idea about the expected additionally travelled distance (up to about 50\% to an exact static mTSP solution).

A continuous approximation model (CAM) for the BD-AVH 
allows to predict the expected travel distance,
when the number of vehicles, number of customers and sequential dynamics scope is given (or can be estimated).
Three models were derived using a Machine Learning approach.
Each offering a different degree in flexibility (ease of use) and accuracy (prediction quality).
Currently, the method was used using uniformly distributed customers in the Euclidean plane.
However, it is possible to use any other stochastic customer location distribution.
The proposed models were derived to operate for medium sized test instances - in a space limited by the number of vehicles, number of customers and scope of dynamics.
These limitations were chosen arbitrarily, 
i.e. setting different limits will return models that can be applied to larger test instances.

This work's algorithms were designed to work with sequential time consideration.
In particular, discrete time events and a finite time horizon were suggested.
However, it would be interesting to consider stochastic time distributions for customers and vehicles.
Focus to the dynamics of vehicles needs to be addressed in future work.
Once these fundamental dynamics are well understood,
Finite State Machines (FSM) and Discrete Event Simulations (DES) 
may proof to be useful tools in developing general dynamic-routing solutions.
Other future investigations could include time-windows for the mTSP
and make used of event-knowledge and schedule timelines.

%% file: tab-set1-abs.tex
 \begin{table*}[htbp]
    \centering
    \caption{BD-AVH and BD-CVH test instances for set 1 using $m$-absolute dynamics.}
 	\begin{tabular}{cc|cccccc}
 		\toprule
 		instance & algorithm & $D_a^m= 0.5$  & $D_a^m=1$    & $D_a^m=1.5$  & $D_a^m=2$    & $D_a^m=4$    & $D_a^m=8$ \\
 		\midrule
 		\multirow{2}[2]{*}{garn9-m2} & BD-AVH & 44.8 & 68.9 & 66.1 & 64.8 & 44.8 & 44.8 \\
 		& BD-CVH & 44.8 & 68.9 & 66.1 & 66.5 & 44.8 & 44.8 \\
 		\midrule
 		\multirow{2}[2]{*}{garn13-m3L4} & BD-AVH & 69.2 & 90.7 & 89.9 & 92.1 & 79.4 & 79.4 \\
 		& BD-CVH & 69.2 & 91.4 & 91   & 93.2 & 87.8 & 87.8 \\
 		\midrule
 		\multirow{2}[2]{*}{garn20-m3} & BD-AVH & 108.1 & 122  & 117.1 & 119.3 & 137.2 & 116.3 \\
 		& BD-CVH & 111.5 & 123.3 & 120  & 120.9 & 140.3 & 120.2 \\
 		\midrule
 		\multirow{2}[2]{*}{bays29-m4} & BD-AVH & 4114 & 4927 & 4131 & 3983 & 3912 & 3458 \\
 		& BD-CVH & 4148 & 5309 & 4260 & 4136 & 3934 & 3643 \\
 		\midrule
 		\multirow{2}[2]{*}{berlin52-m5} & BD-AVH & 22.6k & 25.7k & 18.9k & 17.1k & 17.4k & 13.8k \\
 		& BD-CVH & 23.2k & 26.8k & 22.0k & 18.2k & 17.6k & 13.7k \\
 		\midrule
 		\multirow{2}[2]{*}{eucl-n100m7} & BD-AVH & 38.3k & 42.8k & 37.7k & 31.6k & 29.6k & 26.4k \\
 		& BD-CVH & 37.3k & 44.7k & 37.4k & 30.7k & 31.4k & 27.9k \\
 		\midrule
 		\multirow{2}[2]{*}{lin318-m20} & BD-AVH & 242.0k & 343.3k & 322.5k & 315.9k & 203.6k & 201.2k \\
 		& BD-CVH & 254.9k & 349.2k & 332.0k & 324.0k & 209.7k & 205.4k \\
 		\bottomrule
 	\end{tabular}%
 	\label{tab:set1-abs}%
 \end{table*}%

%% file: tab-set2-rel.tex
\begin{table*}[htbp]
	\centering
	\caption{BD-CVH and BD-AVH set 2 test instances for relative dynamics.}
	\adjustbox{width=\textwidth}{ 
		\begin{tabular}{rcc|ccccccc|cccccccc|}
			\toprule
			\multicolumn{1}{c}{a} & instance & m    & \textbf{2\%} & \textbf{5\%} & \textbf{7\%} & \textbf{10\%} & \textbf{20\%} & \textbf{30\%} & \textbf{100\%} & a    & \textbf{2\%} & \textbf{5\%} & \textbf{7\%} & \textbf{10\%} & \textbf{20\%} & \textbf{30\%} & \textbf{100\%} \\
			\midrule
			\multicolumn{1}{r}{\multirow{24}[12]{*}{\begin{sideways}\textbf{BD-AVH}\end{sideways}}} & eil51 & 2    & 1,251.6 & 1,124.8 & 944.4 & 992.4 & 831.7 & 798.3 & 609.2 & \multirow{24}[12]{*}{\begin{sideways}\textbf{BD-CVH}\end{sideways}} & 1,251.6 & 1,202.9 & 1,069.8 & 1,031.4 & 792.4 & 811.8 & 605.0 \\
			&      & 3    & 1,036.9 & 1,369.6 & 1,112.7 & 957.6 & 935.8 & 892.0 & 656.1 &      & 1,036.9 & 1,507.1 & 1,168.5 & 966.6 & 983.9 & 907.2 & 622.0 \\
			&      & 4    & 1,112.1 & 1,323.0 & 1,224.5 & 1,074.5 & 912.0 & 704.0 & 689.1 &      & 1,112.1 & 1,382.2 & 1,269.7 & 1,123.8 & 1,010.4 & 716.6 & 691.2 \\
			&      & 5    & 1,036.7 & 1,194.4 & 1,192.0 & 1,236.7 & 1,036.4 & 799.6 & 725.4 &      & 1,036.7 & 1,304.5 & 1,212.9 & 1,319.9 & 1,019.1 & 852.4 & 800.4 \\
			\cmidrule{2-10}\cmidrule{12-18}           & eil76 & 2    & 2,295.5 & 1,595.4 & 1,396.9 & 1,275.1 & 960.0 & 938.1 & 766.9 &      & 2,339.8 & 1,519.5 & 1,328.9 & 1,355.9 & 993.3 & 798.1 & 781.2 \\
			&      & 3    & 1,870.5 & 1,650.4 & 1,430.1 & 1,255.4 & 1,025.9 & 970.7 & 725.1 &      & 1,947.4 & 1,677.4 & 1,461.4 & 1,298.5 & 1,127.0 & 1,008.4 & 938.5 \\
			&      & 4    & 1,636.7 & 1,897.6 & 1,623.2 & 1,263.2 & 1,109.7 & 947.6 & 801.0 &      & 1,765.6 & 2,076.3 & 1,725.7 & 1,297.7 & 1,274.5 & 1,077.9 & 845.1 \\
			&      & 5    & 1,754.1 & 1,735.1 & 1,818.6 & 1,388.1 & 1,378.4 & 1,004.3 & 901.4 &      & 1,584.5 & 1,715.4 & 1,958.1 & 1,446.7 & 1,180.3 & 1,037.7 & 948.1 \\
			\cmidrule{2-10}\cmidrule{12-18}           & eil101 & 2    & 2,529.8 & 1,790.7 & 1,666.6 & 1,440.0 & 1,131.1 & 1,043.5 & 822.7 &      & 2,595.2 & 1,846.9 & 1,821.7 & 1,458.4 & 1,186.3 & 1,027.1 & 851.2 \\
			&      & 3    & 2,289.1 & 2,207.0 & 1,877.7 & 1,589.0 & 1,344.4 & 1,020.5 & 881.2 &      & 2,333.3 & 2,248.2 & 1,871.7 & 1,667.7 & 1,220.1 & 1,018.7 & 886.6 \\
			&      & 4    & 2,084.0 & 2,427.6 & 2,027.7 & 1,873.0 & 1,365.2 & 1,102.9 & 990.1 &      & 2,197.6 & 2,552.4 & 2,090.8 & 1,701.7 & 1,389.0 & 1,097.4 & 968.8 \\
			&      & 5    & 1,992.5 & 3,197.5 & 2,235.3 & 1,858.9 & 1,413.4 & 1,172.8 & 1,033.8 &      & 2,081.9 & 3,334.2 & 2,309.1 & 1,945.0 & 1,351.1 & 1,239.1 & 993.6 \\
			\cmidrule{2-10}\cmidrule{12-18}           & kroA100 & 2    & 133.3k & 90.2k & 83.6k & 62.8k & 50.7k & 40.9k & 27.4k &      & 135.3k & 93.7k & 84.0k & 63.2k & 51.2k & 41.5k & 27.4k \\
			&      & 3    & 107.9k & 93.1k & 89.9k & 62.7k & 61.3k & 50.9k & 33.9k &      & 108.9k & 95.4k & 85.2k & 78.0k & 61.4k & 52.9k & 36.1k \\
			&      & 4    & 101.3k & 100.9k & 83.1k & 78.6k & 60.4k & 53.7k & 41.1k &      & 99.6k & 105.0k & 82.1k & 68.8k & 58.6k & 55.6k & 40.3k \\
			&      & 5    & 91.0k & 110.2k & 92.3k & 76.3k & 57.7k & 62.1k & 46.8k &      & 94.5k & 116.4k & 97.9k & 70.5k & 58.7k & 62.1k & 47.6k \\
			\cmidrule{2-10}\cmidrule{12-18}           & kroA150 & 2    & 159.7k & 119.8k & 95.2k & 82.6k & 59.1k & 47.9k & 36.6k &      & 166.7k & 123.1k & 99.2k & 84.1k & 60.8k & 50.8k & 36.7k \\
			&      & 3    & 172.0k & 128.3k & 94.7k & 83.2k & 64.2k & 52.4k & 40.4k &      & 183.3k & 120.7k & 103.5k & 84.8k & 83.0k & 52.8k & 40.4k \\
			&      & 4    & 132.9k & 116.6k & 98.3k & 95.8k & 80.1k & 57.7k & 45.9k &      & 142.5k & 110.0k & 105.1k & 91.0k & 78.7k & 58.4k & 45.3k \\
			&      & 5    & 116.9k & 131.6k & 109.8k & 98.0k & 84.0k & 71.4k & 55.5k &      & 125.6k & 131.4k & 97.0k & 90.7k & 82.9k & 65.3k & 57.3k \\
			\cmidrule{2-10}\cmidrule{12-18}           & kroA200 & 2    & 187.3k & 128.9k & 109.3k & 95.2k & 72.9k & 58.2k & 40.3k &      & 192.8k & 130.6k & 102.2k & 95.2k & 68.4k & 58.4k & 40.4k \\
			&      & 3    & 195.8k & 128.4k & 106.9k & 99.0k & 62.6k & 64.8k & 46.7k &      & 201.7k & 130.8k & 99.6k & 97.4k & 74.7k & 71.1k & 46.9k \\
			&      & 4    & 202.6k & 146.5k & 122.8k & 95.8k & 81.0k & 66.4k & 48.4k &      & 215.3k & 133.0k & 103.1k & 115.4k & 74.4k & 64.6k & 48.7k \\
			&      & 5    & 165.4k & 129.2k & 132.3k & 97.5k & 77.6k & 75.3k & 55.9k &      & 169.0k & 143.8k & 133.9k & 112.4k & 76.1k & 65.1k & 55.6k \\
			\bottomrule
		\end{tabular}%
	}
	\label{tab:set2-rel}%
\end{table*}%

%% file: tab-set2-abs.tex
\begin{table*}[htbp]
	\caption{BD-CVH and BD-AVH set 2 test instances for $m$-absolute dynamics.}
	\adjustbox{width=\textwidth}{ 
		\begin{tabular}{ccccccccc|ccccccccc}
			\toprule
			&      & \boldmath{}\textbf{$D_a^m$}\unboldmath{} & \textbf{0.5} & \textbf{1} & \textbf{1.5} & \textbf{2} & \textbf{4} & \textbf{8} &      &      & \boldmath{}\textbf{$D_a^m$}\unboldmath{} & \textbf{0.5} & \textbf{1} & \textbf{1.5} & \textbf{2} & \textbf{4} & \textbf{8} \\
			\midrule
			& algo. & $m$  & $D_a = 1$ & $D_a = 2$ & $D_a = 3$ & $D_a = 4$ & $D_a = 8$ & $D_a = 16$ &      & algo. & $m$  & $D_a = 1$ & $D_a = 2$ & $D_a = 3$ & $D_a = 4$ & $D_a = 8$ & $D_a = 16$ \\
			\midrule
			\multirow{11}[2]{*}{\begin{sideways}eil51\end{sideways}} & BD-AVH & 2    &   1,251.6  &   1,374.8  &   1,124.8  &      944.4  &      717.7  &      718.6  & \multirow{11}[2]{*}{\begin{sideways}kroA100\end{sideways}} & BD-AVH & 2    & 120.9k & 133.3k & 106.3k & 94.4k & 73.3k & 51.9k \\
			& BD-CVH & 2    &   1,251.6  &   1,436.5  &   1,202.9  &   1,069.8  &      717.7  &      719.5  &      & BD-CVH & 2    & 120.9k & 135.3k & 111.1k & 99.9k & 71.6k & 53.0k \\
			&      & $m$  & $D_a = 2$ & $D_a = 3$ & $D_a = 5$ & $D_a = 6$ & $D_a = 12$ & $D_a = 24$ &      &      & $m$  & $D_a = 2$ & $D_a = 3$ & $D_a = 5$ & $D_a = 6$ & $D_a = 12$ & $D_a = 24$ \\
			& BD-AVH & 3    &   1,185.8  &   1,369.6  &      957.6  &      990.7  &      886.0  &      650.9  &      & BD-AVH & 3    & 107.9k & 111.8k & 93.1k & 83.2k & 67.0k & 55.0k \\
			& BD-CVH & 3    &   1,236.3  &   1,507.1  &      966.6  &   1,007.7  &      891.2  &      664.0  &      & BD-CVH & 3    & 108.9k & 119.9k & 95.4k & 85.3k & 69.3k & 54.4k \\
			&      & $m$  & $D_a = 2$ & $D_a = 4$ & $D_a = 6$ & $D_a = 8$ & $D_a = 16$ & $D_a = 32$ &      &      & $m$  & $D_a = 2$ & $D_a = 4$ & $D_a = 6$ & $D_a = 8$ & $D_a = 16$ & $D_a = 32$ \\
			& BD-AVH & 4    &   1,143.4  &   1,224.5  &   1,057.4  &   1,002.6  &      680.6  &      784.6  &      & BD-AVH & 4    & 101.3k & 101.3k & 85.9k & 79.2k & 63.9k & 55.8k \\
			& BD-CVH & 4    &   1,116.9  &   1,269.7  &   1,075.9  &   1,064.9  &      661.0  &      705.6  &      & BD-CVH & 4    & 99.6k & 107.1k & 90.7k & 80.7k & 62.3k & 52.8k \\
			&      & $m$  & $D_a = 3$ & $D_a = 5$ & $D_a = 8$ & $D_a = 10$ & $D_a = 20$ & $D_a = 40$ &      &      & $m$  & $D_a = 3$ & $D_a = 5$ & $D_a = 8$ & $D_a = 10$ & $D_a = 20$ & $D_a = 40$ \\
			& BD-AVH & 5    &   1,194.4  &   1,236.7  &   1,117.5  &   1,036.4  &      691.3  &      725.5  &      & BD-AVH & 5    & 83.4k & 110.2k & 87.7k & 76.3k & 57.7k & 57.4k \\
			& BD-CVH & 5    &   1,304.5  &   1,319.9  &   1,121.7  &   1,019.1  &      698.7  &      748.3  &      & BD-CVH & 5    & 88.8k & 116.4k & 82.1k & 70.5k & 58.7k & 57.7k \\
			\midrule
			& algo. & $m$  & $D_a = 1$ & $D_a = 2$ & $D_a = 3$ & $D_a = 4$ & $D_a = 8$ & $D_a = 16$ &      & algo. & $m$  & $D_a = 1$ & $D_a = 2$ & $D_a = 3$ & $D_a = 4$ & $D_a = 8$ & $D_a = 16$ \\
			\midrule
			\multirow{11}[2]{*}{\begin{sideways}eil76\end{sideways}} & BD-AVH & 2    &   1,727.1  &   2,295.5  &   1,737.0  &   1,595.4  &   1,275.1  &   1,048.4  & \multirow{11}[2]{*}{\begin{sideways}kroA150\end{sideways}} & BD-AVH & 2    & 170.5k & 205.8k & 159.7k & 142.8k & 105.3k & 78.4k \\
			& BD-CVH & 2    &   1,727.1  &   2,339.8  &   1,755.9  &   1,519.5  &   1,355.9  &      987.0  &      & BD-CVH & 2    & 170.5k & 208.6k & 166.7k & 143.4k & 104.0k & 77.7k \\
			&      & $m$  & $D_a = 2$ & $D_a = 3$ & $D_a = 5$ & $D_a = 6$ & $D_a = 12$ & $D_a = 24$ &      &      & $m$  & $D_a = 2$ & $D_a = 3$ & $D_a = 5$ & $D_a = 6$ & $D_a = 12$ & $D_a = 24$ \\
			& BD-AVH & 3    &   1,870.5  &   2,221.4  &   1,430.1  &   1,483.1  &   1,125.6  &      921.7  &      & BD-AVH & 3    & 157.7k & 172.0k & 125.8k & 120.0k & 83.9k & 80.5k \\
			& BD-CVH & 3    &   1,947.4  &   2,332.3  &   1,461.4  &   1,418.9  &   1,006.8  &      845.3  &      & BD-CVH & 3    & 159.3k & 183.3k & 133.4k & 125.7k & 88.2k & 82.6k \\
			&      & $m$  & $D_a = 2$ & $D_a = 4$ & $D_a = 6$ & $D_a = 8$ & $D_a = 16$ & $D_a = 32$ &      &      & $m$  & $D_a = 2$ & $D_a = 4$ & $D_a = 6$ & $D_a = 8$ & $D_a = 16$ & $D_a = 32$ \\
			& BD-AVH & 4    &   1,636.7  &   1,897.6  &   1,372.3  &   1,263.2  &   1,075.0  &      882.6  &      & BD-AVH & 4    & 134.8k & 152.3k & 119.7k & 113.6k & 93.5k & 75.6k \\
			& BD-CVH & 4    &   1,765.6  &   2,076.3  &   1,395.7  &   1,297.7  &   1,224.2  &      877.5  &      & BD-CVH & 4    & 137.0k & 160.5k & 126.1k & 115.5k & 88.5k & 61.5k \\
			&      & $m$  & $D_a = 3$ & $D_a = 5$ & $D_a = 8$ & $D_a = 10$ & $D_a = 20$ & $D_a = 40$ &      &      & $m$  & $D_a = 3$ & $D_a = 5$ & $D_a = 8$ & $D_a = 10$ & $D_a = 20$ & $D_a = 40$ \\
			& BD-AVH & 5    &   1,799.8  &   1,818.6  &   1,388.1  &   1,434.4  &   1,134.0  &      873.8  &      & BD-AVH & 5    & 116.9k & 157.3k & 121.4k & 109.8k & 80.5k & 66.1k \\
			& BD-CVH & 5    &   1,828.2  &   1,958.1  &   1,446.7  &   1,294.7  &   1,128.0  &      848.7  &      & BD-CVH & 5    & 125.6k & 167.9k & 113.9k & 97.0k & 81.6k & 71.0k \\
			\midrule
			& algo. & $m$  & $D_a = 1$ & $D_a = 2$ & $D_a = 3$ & $D_a = 4$ & $D_a = 8$ & $D_a = 16$ &      & algo. & $m$  & $D_a = 1$ & $D_a = 2$ & $D_a = 3$ & $D_a = 4$ & $D_a = 8$ & $D_a = 16$ \\
			\midrule
			\multirow{11}[2]{*}{\begin{sideways}eil101\end{sideways}} & BD-AVH & 2    &   1,836.8  &   2,529.8  &   2,278.5  &   2,043.3  &   1,551.4  &   1,218.8  & \multirow{11}[2]{*}{\begin{sideways}kroA200\end{sideways}} & BD-AVH & 2    & 220.9k & 287.7k & 208.8k & 187.3k & 150.9k & 99.9k \\
			& BD-CVH & 2    &   1,836.8  &   2,595.2  &   2,361.8  &   2,073.8  &   1,625.4  &   1,187.5  &      & BD-CVH & 2    & 220.9k & 291.8k & 215.0k & 192.8k & 147.0k & 106.4k \\
			&      & $m$  & $D_a = 2$ & $D_a = 3$ & $D_a = 5$ & $D_a = 6$ & $D_a = 12$ & $D_a = 24$ &      &      & $m$  & $D_a = 2$ & $D_a = 3$ & $D_a = 5$ & $D_a = 6$ & $D_a = 12$ & $D_a = 24$ \\
			& BD-AVH & 3    &   2,289.1  &   2,627.0  &   2,207.0  &   1,911.3  &   1,527.3  &   1,043.3  &      & BD-AVH & 3    & 207.3k & 241.6k & 174.7k & 155.2k & 119.1k & 78.1k \\
			& BD-CVH & 3    &   2,333.3  &   2,722.8  &   2,248.2  &   1,908.5  &   1,478.4  &   1,170.1  &      & BD-CVH & 3    & 211.7k & 249.1k & 180.4k & 153.9k & 113.3k & 78.7k \\
			&      & $m$  & $D_a = 2$ & $D_a = 4$ & $D_a = 6$ & $D_a = 8$ & $D_a = 16$ & $D_a = 32$ &      &      & $m$  & $D_a = 2$ & $D_a = 4$ & $D_a = 6$ & $D_a = 8$ & $D_a = 16$ & $D_a = 32$ \\
			& BD-AVH & 4    &   2,084.0  &   2,851.2  &   2,099.5  &   1,868.2  &   1,472.1  &   1,162.8  &      & BD-AVH & 4    & 175.1k & 202.6k & 168.3k & 139.8k & 105.1k & 76.2k \\
			& BD-CVH & 4    &   2,197.6  &   2,988.7  &   2,402.3  &   1,867.7  &   1,648.3  &   1,129.6  &      & BD-CVH & 4    & 176.2k & 215.3k & 166.1k & 140.1k & 105.9k & 75.4k \\
			&      & $m$  & $D_a = 3$ & $D_a = 5$ & $D_a = 8$ & $D_a = 10$ & $D_a = 20$ & $D_a = 40$ &      &      & $m$  & $D_a = 3$ & $D_a = 5$ & $D_a = 8$ & $D_a = 10$ & $D_a = 20$ & $D_a = 40$ \\
			& BD-AVH & 5    &   2,167.3  &   3,197.5  &   2,078.5  &   1,858.9  &   1,413.4  &   1,117.0  &      & BD-AVH & 5    & 178.7k & 204.8k & 158.2k & 129.2k & 97.5k & 77.6k \\
			& BD-CVH & 5    &   2,381.8  &   3,334.2  &   2,180.9  &   1,945.0  &   1,351.1  &   1,187.2  &      & BD-CVH & 5    & 180.7k & 217.1k & 159.6k & 143.8k & 112.4k & 76.1k \\
			\bottomrule
		\end{tabular}%
	}
	\label{tab:set2-abs}%
\end{table*}%

%% file: tab-setX-rel.tex
\begin{table*}[htbp]
	\centering
	\caption{BD-CVH and BD-AVH set X-U test instances for relative dynamics.}
	\adjustbox{width=\textwidth}{ 	
		\begin{tabular}{cc|c|ccccccc|cccccccc}
		\toprule
		a    & instance & \multicolumn{1}{c|}{m} & \textbf{$D_r=2\%$} & \textbf{$D_r=5\%$} & \textbf{$D_r=7\%$} & \textbf{$D_r=10\%$} & \textbf{$D_r=20\%$} & \textbf{$D_r=30\%$} & \textbf{$D_r=100\%$} & a    & \textbf{$D_r=2\%$} & \textbf{$D_r=5\%$} & \textbf{$D_r=7\%$} & \textbf{$D_r=10\%$} & \textbf{$D_r=20\%$} & \textbf{$D_r=30\%$} & \multicolumn{1}{c}{\textbf{$D_r=100\%$}} \\
		\midrule
		\multirow{16}[4]{*}{\begin{sideways}BD-AVH\end{sideways}} & X-n115-k10 & 10   & 25.3k & 28.6k & 30.9k & 31.3k & 24.7k & 20.8k & 18.9k & \multirow{16}[4]{*}{\begin{sideways}BD-CVH\end{sideways}} & 26.0k & 29.5k & 32.4k & 34.4k & 25.3k & 20.5k & 19.4k \\
		& X-n153-k22 & 22   & 127.3M & 78.4M & 39.2M & 19.6M & 30.1k & 34.9k & 38.7k &      & 117.5M & 68.6M & 39.2M & 19.6M & 29.7k & 37.6k & 38.6k \\
		& X-n176-k26 & 26   & 240.3M & 160.2M & 112.2M & 32.1M & 55.5k & 58.5k & 62.4k &      & 240.3M & 160.2M & 128.2M & 16.1M & 59.2k & 61.1k & 63.0k \\
		& X-n214-k11 & 11   & 23.8k & 27.8k & 24.8k & 21.4k & 18.1k & 15.3k & 16.2k &      & 24.1k & 31.2k & 27.5k & 22.2k & 18.1k & 17.0k & 16.5k \\
		& X-n233-k16 & 16   & 37.3k & 44.4k & 53.0k & 39.5k & 35.6k & 39.1k & 33.1k &      & 38.2k & 46.9k & 58.9k & 42.4k & 32.7k & 37.0k & 35.4k \\
		& X-n270-k35 & 35   & 572.9M & 191.0M & 38.2M & 0.0M & 46.3k & 49.4k & 63.4k &      & 572.9M & 267.4M & 76.4M & 54.0k & 48.6k & 48.9k & 65.3k \\
		& X-n313-k71 & 71   & 2.51G & 2.16G & 1.92G & 1.38G & 106.3k & 110.4k & 139.8k &      & 2.55G & 2.11G & 1.82G & 1.38G & 110.9k & 114.9k & 141.8k \\
		& X-n327-k20 & 20   & 56.2k & 61.6k & 65.3k & 56.0k & 45.5k & 42.4k & 47.1k &      & 59.3k & 64.8k & 71.0k & 56.2k & 47.3k & 43.6k & 44.0k \\
		& X-n367-k17 & 17   & 157.3M & 61.3k & 49.1k & 41.7k & 44.7k & 41.8k & 35.0k &      & 157.3M & 66.1k & 55.0k & 44.8k & 43.1k & 42.5k & 35.8k \\
		\cmidrule{2-10}\cmidrule{12-18}           & X-n429-k61 & 61   & 3.41G & 2.34G & 1.27G & 194.8M & 93.3k & 121.8k & 116.8k &      & 3.21G & 2.24G & 1.27G & 194.8M & 97.5k & 129.1k & 118.7k \\
		& X-n491-k59 & 59   & 3.53G & 1.22G & 86.7k & 94.0k & 94.5k & 90.6k & 105.1k &      & 3.41G & 1,216.2M & 90.9k & 100.3k & 105.7k & 92.7k & 107.0k \\
		& X-n536-k96 & 96   & 11.11G & 7.95G & 6.75G & 4.50G & 112.5k & 127.6k & 185.4k &      & 10.8G & 7.80G & 6.60G & 4.35G & 119.5k & 131.9k & 185.6k \\
		& X-n641-k35 & 35   & 2.48G & 125.0k & 113.1k & 110.1k & 100.8k & 105.2k & 93.6k &      & 2.27G & 133.0k & 124.8k & 112.2k & 106.4k & 107.4k & 98.7k \\
		& X-n783-k48 & 48   & 3.55G & 135.4k & 138.7k & 125.2k & 107.7k & 119.1k & 118.0k &      & 3.87G & 142.2k & 151.3k & 133.8k & 117.2k & 118.4k & 121.0k \\
		& X-n837-k142 & 142  & 37.8G & 29.4G & 24.1G & 14.3G & 238.1k & 260.5k & 294.9k &      & 37.8G & 29.4G & 23.4G & 12.9G & 253.1k & 272.9k & 300.5k \\
		& X-n936-k151 & 151  & 45.3G & 32.6G & 27.2G & 13.1G & 156.1k & 176.9k & 181.5k &      & 45.3G & 31.7G & 26.3G & 14.5G & 165.7k & 192.8k & 187.3k \\
		\bottomrule
	\end{tabular}%
}
	\label{tab:setX-rel}%
\end{table*}%

%% file: tab-setX-abs.tex
\begin{table*}[htbp]
	\centering
	\caption{BD-AVH and BD-CVH test instances for set X-U using $m$-absolute dynamics.}
	\adjustbox{totalheight=.9\textheight-2\baselineskip}{ 	
	\begin{tabular}{cccccccc}
		\toprule
		& \boldmath{}\textbf{$D_a^m$}\unboldmath{} & \textbf{0.5} & \textbf{1} & \textbf{1.5} & \textbf{2} & \textbf{4} & \multicolumn{1}{c}{\textbf{8}} \\
		\midrule
		X-n115-k10 & $m$  & $D_a = 5$ & $D_a = 10$ & $D_a = 15$ & $D_a = 20$ & $D_a = 40$ & $D_a = 80$ \\
		\midrule
		BD-AVH & 10   & 27.3k & 34.6k & 25.2k & 23.8k & 17.9k & 18.5k \\
		BD-CVH & 10   & 27.9k & 37.3k & 26.8k & 24.4k & 22.5k & 18.6k \\
		\midrule
		X-n153-k22 & $m$  & $D_a = 11$ & $D_a = 22$ & $D_a = 33$ & $D_a = 44$ & $D_a = 88$ & $D_a = 176$ \\
		\midrule
		BD-AVH & 22   & 39.2M & 35.4k & 30.0k & 33.6k & 37.2k & 38.7k \\
		BD-CVH & 22   & 39.2M & 37.6k & 30.3k & 38.4k & 38.1k & 38.6k \\
		\midrule
		X-n176-k26 & $m$  & $D_a = 13$ & $D_a = 26$ & $D_a = 39$ & $D_a = 52$ & $D_a = 104$ & $D_a = 208$ \\
		\midrule
		BD-AVH & 26   & 80.1M & 69.2k & 59.1k & 59.0k & 64.8k & 62.4k \\
		BD-CVH & 26   & 96.2M & 72.6k & 58.9k & 59.5k & 65.6k & 63.0k \\
		\midrule
		X-n214-k11 & $m$  & $D_a = 6$ & $D_a = 11$ & $D_a = 17$ & $D_a = 22$ & $D_a = 44$ & $D_a = 88$ \\
		\midrule
		BD-AVH & 11   & 24.6k & 27.8k & 22.0k & 21.3k & 17.9k & 15.3k \\
		BD-CVH & 11   & 25.1k & 31.2k & 24.0k & 22.0k & 19.4k & 15.6k \\
		\midrule
		X-n233-k16 & $m$  & $D_a = 8$ & $D_a = 16$ & $D_a = 24$ & $D_a = 32$ & $D_a = 64$ & $D_a = 128$ \\
		\midrule
		BD-AVH & 16   & 37.9k & 53.0k & 42.4k & 34.0k & 31.7k & 33.0k \\
		BD-CVH & 16   & 38.5k & 58.9k & 44.4k & 36.4k & 33.6k & 33.3k \\
		\midrule
		X-n270-k35 & $m$  & $D_a = 18$ & $D_a = 35$ & $D_a = 53$ & $D_a = 70$ & $D_a = 140$ & $D_a = 280$ \\
		\midrule
		BD-AVH & 35   & 46.9k & 59.0k & 46.4k & 51.9k & 49.4k & 63.4k \\
		BD-CVH & 35   & 47.6k & 64.0k & 50.5k & 50.5k & 50.9k & 65.3k \\
		\midrule
		X-n313-k71 & $m$  & $D_a = 36$ & $D_a = 71$ & $D_a = 107$ & $D_a = 142$ & $D_a = 284$ & $D_a = 568$ \\
		\midrule
		BD-AVH & 71   & 1,130.1M & 119.7k & 110.5k & 122.4k & 138.0k & 139.8k \\
		BD-CVH & 71   & 1,130.1M & 123.9k & 118.0k & 126.2k & 141.8k & 141.8k \\
		\midrule
		X-n327-k20 & $m$  & $D_a = 10$ & $D_a = 20$ & $D_a = 30$ & $D_a = 40$ & $D_a = 80$ & $D_a = 160$ \\
		\midrule
		BD-AVH & 20   & 57.2k & 74.5k & 54.9k & 52.0k & 42.9k & 49.7k \\
		BD-CVH & 20   & 60.9k & 81.7k & 60.7k & 53.6k & 43.2k & 51.0k \\
		\midrule
		X-n367-k17 & $m$  & $D_a = 9$ & $D_a = 17$ & $D_a = 26$ & $D_a = 34$ & $D_a = 68$ & $D_a = 136$ \\
		\midrule
		BD-AVH & 17   & 47.4k & 64.2k & 49.1k & 46.5k & 44.8k & 41.6k \\
		BD-CVH & 17   & 49.9k & 70.3k & 55.0k & 43.9k & 41.1k & 41.6k \\
		\midrule
		X-n429-k61 & $m$  & $D_a = 31$ & $D_a = 61$ & $D_a = 92$ & $D_a = 122$ & $D_a = 244$ & $D_a = 488$ \\
		\midrule
		BD-AVH & 61   & 1,265.9M & 103.7k & 93.6k & 113.6k & 115.2k & 116.8k \\
		BD-CVH & 61   & 1,265.9M & 112.8k & 98.2k & 117.2k & 122.1k & 118.7k \\
		\midrule
		X-n491-k59 & $m$  & $D_a = 30$ & $D_a = 59$ & $D_a = 89$ & $D_a = 118$ & $D_a = 236$ & $D_a = 472$ \\
		\midrule
		BD-AVH & 59   & 486.5M & 107.4k & 97.5k & 102.5k & 109.9k & 108.3k \\
		BD-CVH & 59   & 729.8M & 116.4k & 104.6k & 107.4k & 113.2k & 106.7k \\
		\midrule
		X-n536-k96 & $m$  & $D_a = 48$ & $D_a = 96$ & $D_a = 144$ & $D_a = 192$ & $D_a = 384$ & $D_a = 768$ \\
		\midrule
		BD-AVH & 96   & 4,502.6M & 121.4k & 113.3k & 138.4k & 183.3k & 185.4k \\
		BD-CVH & 96   & 4,352.5M & 127.2k & 123.9k & 141.6k & 186.3k & 185.6k \\
		\midrule
		X-n641-k35 & $m$  & $D_a = 18$ & $D_a = 35$ & $D_a = 53$ & $D_a = 70$ & $D_a = 140$ & $D_a = 280$ \\
		\midrule
		BD-AVH & 35   & 1,240.5M & 132.3k & 109.1k & 108.0k & 100.3k & 97.1k \\
		BD-CVH & 35   & 1,240.5M & 146.7k & 116.9k & 116.4k & 110.9k & 102.1k \\
		\midrule
		X-n783-k48 & $m$  & $D_a = 24$ & $D_a = 48$ & $D_a = 72$ & $D_a = 96$ & $D_a = 192$ & $D_a = 384$ \\
		\midrule
		BD-AVH & 48   & 1,290.4M & 170.5k & 135.7k & 116.4k & 121.1k & 130.6k \\
		BD-CVH & 48   & 967.8M & 190.8k & 144.2k & 124.9k & 121.7k & 134.8k \\
		\midrule
		X-n837-k142 & $m$  & $D_a = 71$ & $D_a = 142$ & $D_a = 213$ & $D_a = 284$ & $D_a = 568$ & $D_a = 1136$ \\
		\midrule
		BD-AVH & 142  & 18.9G & 247.3k & 244.2k & 275.9k & 298.3k & 294.9k \\
		BD-CVH & 142  & 17.8G & 261.5k & 256.2k & 282.2k & 300.4k & 300.5k \\
		\midrule
		X-n936-k151 & $m$  & $D_a = 76$ & $D_a = 151$ & $D_a = 227$ & $D_a = 302$ & $D_a = 604$ & $D_a = 1208$ \\
		\midrule
		BD-AVH & 151  & 20.8G & 178.8k & 162.6k & 175.8k & 187.4k & 181.5k \\
		BD-CVH & 151  & 22.2G & 199.2k & 180.2k & 187.3k & 194.1k & 187.3k \\
		\bottomrule
	\end{tabular}%
}
	\label{tab:setX-abs}%
\end{table*}%

%% file: tab-setXXL-abs.tex
\begin{table*}[htbp]
	\centering
	\caption{BD-CVH test instances for set XXL using $m$-absolute dynamics}
	\begin{tabular}{ccccc|cccccc}
		\toprule
		\textbf{ instance} & \textbf{n} & \textbf{m} & \textbf{L} & \textbf{attribute} & \boldmath{}\textbf{$D_a^m = 8$}\unboldmath{} & \boldmath{}\textbf{$D_a^m = 4$}\unboldmath{} & \boldmath{}\textbf{$D_a^m = 2$}\unboldmath{} & \boldmath{}\textbf{$D_a^m = 1.5$}\unboldmath{} & \boldmath{}\textbf{$D_a^m = 1$}\unboldmath{} & \boldmath{}\textbf{$D_a^m = 0.5$}\unboldmath{} \\
		\midrule
		\multirow{3}[2]{*}{Flanders2} & \multirow{3}[2]{*}{30k} & \multirow{3}[2]{*}{268} & \multirow{3}[2]{*}{112} & $D_a$ & 2,144 & 1,072 & 536  & 402  & 268  & 134 \\
		&      &      &      & distance & 18.7M & 13.9M & 15.6M & 18.4M & 35.2M & 23.3M \\
		&      &      &      & seconds & 34.0 & 19.0 & 10.8 & 9.3  & 7.6  & 5.0 \\
		\midrule
		\multirow{3}[2]{*}{Flanders1} & \multirow{3}[2]{*}{20k} & \multirow{3}[2]{*}{717} & \multirow{3}[2]{*}{28} & $D_a$ & 5,736 & 2,868 & 1,434 & 1,076 & 717  & 359 \\
		&      &      &      & distance & 26.5M & 28.2M & 19.7M & 18.7M & 19.7M & 16.6M \\
		&      &      &      & seconds & 459.7 & 210.2 & 45.1 & 32.1 & 21.4 & 9.1 \\
		\midrule
		\multirow{3}[2]{*}{Brussels2} & \multirow{3}[2]{*}{16k} & \multirow{3}[2]{*}{190} & \multirow{3}[2]{*}{85} & $D_a$ & 1,520 & 760  & 380  & 285  & 190  & 95 \\
		&      &      &      & distance & 819k & 1.02M & 1.18M & 1.31M & 2.25M & 1.72M \\
		&      &      &      & seconds & 10.2 & 5.4  & 3.0  & 2.5  & 1.9  & 1.4 \\
		\midrule
		\multirow{3}[2]{*}{Brussels1} & \multirow{3}[2]{*}{15k} & \multirow{3}[2]{*}{537} & \multirow{3}[2]{*}{28} & $D_a$ & 4,296 & 2,148 & 1,074 & 806  & 537  & 269 \\
		&      &      &      & distance & 1.37M & 1.30M & 1.18M & 1.13M & 1.61M & 1.37M \\
		&      &      &      & seconds & 173.8 & 52.8 & 17.0 & 13.1 & 9.0  & 3.9 \\
		\midrule
		\multirow{3}[2]{*}{Ghent2} & \multirow{3}[2]{*}{11k} & \multirow{3}[2]{*}{115} & \multirow{3}[2]{*}{96} & $D_a$ & 920  & 460  & 230  & 173  & 115  & 58 \\
		&      &      &      & distance & 599k & 685k & 848k & 978k & 1.72M & 1.34M \\
		&      &      &      & seconds & 2.7  & 1.6  & 1.0  & 0.8  & 0.8  & 0.6 \\
		\midrule
		\multirow{3}[2]{*}{Ghent1} & \multirow{3}[2]{*}{10k} & \multirow{3}[2]{*}{509} & \multirow{3}[2]{*}{20} & $D_a$ & 4072 & 2036 & 1018 & 764  & 509  & 255 \\
		&      &      &      & distance & 1.20M & 1.20M & 1.08M & 846k & 1.17M & 1.04M \\
		&      &      &      & seconds & 89.2 & 30.2 & 10.2 & 7.7  & 5.3  & 2.3 \\
		\midrule
		\multirow{3}[2]{*}{Antwerp2} & \multirow{3}[2]{*}{7k} & \multirow{3}[2]{*}{125} & \multirow{3}[2]{*}{57} & $D_a$ & 1,000 & 500  & 250  & 188  & 125  & 63 \\
		&      &      &      & distance & 778k & 768k & 873k & 978k & 1.61M & 1.15M \\
		&      &      &      & seconds & 1.7  & 1.0  & 0.6  & 0.5  & 0.4  & 0.3 \\
		\midrule
		\multirow{3}[2]{*}{Antwerp1} & \multirow{3}[2]{*}{6k} & \multirow{3}[2]{*}{359} & \multirow{3}[2]{*}{17} & $D_a$ & 2,872 & 1,436 & 718  & 539  & 359  & 180 \\
		&      &      &      & distance & 1.13M & 1.13M & 1.12M & 0.97M & 1.19M & 1.02M \\
		&      &      &      & seconds & 13.2 & 5.96 & 3.03 & 2.32 & 1.58 & 0.82 \\
		\midrule
		\multirow{3}[2]{*}{Leuven2} & \multirow{3}[2]{*}{4k} & \multirow{3}[2]{*}{47} & \multirow{3}[2]{*}{86} & $D_a$ & 376  & 188  & 94   & 71   & 47   & 24 \\
		&      &      &      & distance & 304k & 307k & 433k & 504k & 787k & 618k \\
		&      &      &      & seconds & 0.22 & 0.16 & 0.11 & 0.11 & 0.10 & 0.10 \\
		\midrule
		\multirow{3}[2]{*}{Leuven1} & \multirow{3}[2]{*}{3k} & \multirow{3}[2]{*}{212} & \multirow{3}[2]{*}{15} & $D_a$ & 1,696 & 848  & 424  & 318  & 212  & 106 \\
		&      &      &      & distance & 515k & 557k & 501k & 404k & 451k & 427k \\
		&      &      &      & seconds & 1.48 & 0.89 & 0.48 & 0.37 & 0.27 & 0.16 \\
		\bottomrule
	\end{tabular}%
	\label{tab:setXXL-abs}%
\end{table*}%

%% file: BD-mTSP.bbl
\begin{thebibliography}{40}
\expandafter\ifx\csname natexlab\endcsname\relax\def\natexlab#1{#1}\fi
\providecommand{\url}[1]{\texttt{#1}}
\providecommand{\href}[2]{#2}
\providecommand{\path}[1]{#1}
\providecommand{\DOIprefix}{doi:}
\providecommand{\ArXivprefix}{arXiv:}
\providecommand{\URLprefix}{URL: }
\providecommand{\Pubmedprefix}{pmid:}
\providecommand{\doi}[1]{\href{http://dx.doi.org/#1}{\path{#1}}}
\providecommand{\Pubmed}[1]{\href{pmid:#1}{\path{#1}}}
\providecommand{\bibinfo}[2]{#2}
\ifx\xfnm\relax \def\xfnm[#1]{\unskip,\space#1}\fi
\bibitem[{Ansari et~al.(2018)Ansari, Ba{\c{s}}dere, Li, Ouyang and
  Smilowitz}]{ansari2018advancements}
\bibinfo{author}{Ansari, S.}, \bibinfo{author}{Ba{\c{s}}dere, M.},
  \bibinfo{author}{Li, X.}, \bibinfo{author}{Ouyang, Y.},
  \bibinfo{author}{Smilowitz, K.}, \bibinfo{year}{2018}.
\newblock \bibinfo{title}{Advancements in continuous approximation models for
  logistics and transportation systems: 1996--2016}.
\newblock \bibinfo{journal}{Transportation Research Part B: Methodological}
  \bibinfo{volume}{107}, \bibinfo{pages}{229--252}.
\bibitem[{Applegate et~al.(2006)Applegate, Bixby, Chvatal and
  Cook}]{applegate2006traveling}
\bibinfo{author}{Applegate, D.L.}, \bibinfo{author}{Bixby, R.E.},
  \bibinfo{author}{Chvatal, V.}, \bibinfo{author}{Cook, W.J.},
  \bibinfo{year}{2006}.
\newblock \bibinfo{title}{The traveling salesman problem: a computational
  study}.
\newblock \bibinfo{publisher}{Princeton university press}.
\bibitem[{Arnold et~al.(2019)Arnold, Gendreau and
  S{\"o}rensen}]{arnold2019efficiently}
\bibinfo{author}{Arnold, F.}, \bibinfo{author}{Gendreau, M.},
  \bibinfo{author}{S{\"o}rensen, K.}, \bibinfo{year}{2019}.
\newblock \bibinfo{title}{Efficiently solving very large-scale routing
  problems}.
\newblock \bibinfo{journal}{Computers \& Operations Research}
  \bibinfo{volume}{107}, \bibinfo{pages}{32--42}.
\bibitem[{Beardwood et~al.(1959)Beardwood, Halton and
  Hammersley}]{beardwood1959shortest}
\bibinfo{author}{Beardwood, J.}, \bibinfo{author}{Halton, J.H.},
  \bibinfo{author}{Hammersley, J.M.}, \bibinfo{year}{1959}.
\newblock \bibinfo{title}{The shortest path through many points}, in:
  \bibinfo{booktitle}{Mathematical Proceedings of the Cambridge Philosophical
  Society}, \bibinfo{organization}{Cambridge University Press}. pp.
  \bibinfo{pages}{299--327}.
\bibitem[{Bekta{\c{s}} et~al.(2019)Bekta{\c{s}}, Gouveia, Mart{\'\i}nez-Sykora
  and Salazar-Gonz{\'a}lez}]{bektacs2019balanced}
\bibinfo{author}{Bekta{\c{s}}, T.}, \bibinfo{author}{Gouveia, L.},
  \bibinfo{author}{Mart{\'\i}nez-Sykora, A.},
  \bibinfo{author}{Salazar-Gonz{\'a}lez, J.J.}, \bibinfo{year}{2019}.
\newblock \bibinfo{title}{Balanced vehicle routing: Polyhedral analysis and
  branch-and-cut algorithm}.
\newblock \bibinfo{journal}{European Journal of Operational Research}
  \bibinfo{volume}{273}, \bibinfo{pages}{452--463}.
\bibitem[{Caramia et~al.(2002)Caramia, Italiano, Oriolo, Pacifici and
  Perugia}]{caramia2002routing}
\bibinfo{author}{Caramia, M.}, \bibinfo{author}{Italiano, G.F.},
  \bibinfo{author}{Oriolo, G.}, \bibinfo{author}{Pacifici, A.},
  \bibinfo{author}{Perugia, A.}, \bibinfo{year}{2002}.
\newblock \bibinfo{title}{Routing a fleet of vehicles for dynamic combined
  pick-up and deliveries services}, in: \bibinfo{booktitle}{Operations Research
  Proceedings 2001}, \bibinfo{organization}{Springer}. pp.
  \bibinfo{pages}{3--8}.
\bibitem[{{\c{C}}avdar and Sokol(2015)}]{ccavdar2015distribution}
\bibinfo{author}{{\c{C}}avdar, B.}, \bibinfo{author}{Sokol, J.},
  \bibinfo{year}{2015}.
\newblock \bibinfo{title}{A distribution-free tsp tour length estimation model
  for random graphs}.
\newblock \bibinfo{journal}{European Journal of Operational Research}
  \bibinfo{volume}{243}, \bibinfo{pages}{588--598}.
\bibitem[{Cheung et~al.(2008)Cheung, Choy, Li, Shi and
  Tang}]{cheung2008dynamic}
\bibinfo{author}{Cheung, B.K.S.}, \bibinfo{author}{Choy, K.},
  \bibinfo{author}{Li, C.L.}, \bibinfo{author}{Shi, W.}, \bibinfo{author}{Tang,
  J.}, \bibinfo{year}{2008}.
\newblock \bibinfo{title}{Dynamic routing model and solution methods for fleet
  management with mobile technologies}.
\newblock \bibinfo{journal}{International Journal of Production Economics}
  \bibinfo{volume}{113}, \bibinfo{pages}{694--705}.
\bibitem[{Christofides and Eilon(1972)}]{christofides1972algorithms}
\bibinfo{author}{Christofides, N.}, \bibinfo{author}{Eilon, S.},
  \bibinfo{year}{1972}.
\newblock \bibinfo{title}{Algorithms for large-scale travelling salesman
  problems}.
\newblock \bibinfo{journal}{Journal of the Operational Research Society}
  \bibinfo{volume}{23}, \bibinfo{pages}{511--518}.
\bibitem[{Cordeau and Laporte(2003)}]{cordeau2003dial}
\bibinfo{author}{Cordeau, J.F.}, \bibinfo{author}{Laporte, G.},
  \bibinfo{year}{2003}.
\newblock \bibinfo{title}{The dial-a-ride problem (darp): Variants, modeling
  issues and algorithms}.
\newblock \bibinfo{journal}{Quarterly Journal of the Belgian, French and
  Italian Operations Research Societies} \bibinfo{volume}{1},
  \bibinfo{pages}{89--101}.
\bibitem[{Daganzo(1984)}]{daganzo1984distance}
\bibinfo{author}{Daganzo, C.F.}, \bibinfo{year}{1984}.
\newblock \bibinfo{title}{The distance traveled to visit n points with a
  maximum of c stops per vehicle: An analytic model and an application}.
\newblock \bibinfo{journal}{Transportation science} \bibinfo{volume}{18},
  \bibinfo{pages}{331--350}.
\bibitem[{Eilon and Christofides(1971)}]{Eilon1971}
\bibinfo{author}{Eilon, Samuel, C.D.T.W.G.}, \bibinfo{author}{Christofides,
  N.}, \bibinfo{year}{1971}.
\newblock \bibinfo{title}{Distribution Management: Mathematical Modelling and
  Practical Analysis}.
\newblock \bibinfo{publisher}{New York, Hafner}.
\bibitem[{Erera and Daganzo(2003)}]{erera2003dynamic}
\bibinfo{author}{Erera, A.L.}, \bibinfo{author}{Daganzo, C.F.},
  \bibinfo{year}{2003}.
\newblock \bibinfo{title}{A dynamic scheme for stochastic vehicle routing}.
\newblock \bibinfo{journal}{Report, Georgia Institute of Technology} .
\bibitem[{Fabri and Recht(2006)}]{fabri2006dynamic}
\bibinfo{author}{Fabri, A.}, \bibinfo{author}{Recht, P.}, \bibinfo{year}{2006}.
\newblock \bibinfo{title}{On dynamic pickup and delivery vehicle routing with
  several time windows and waiting times}.
\newblock \bibinfo{journal}{Transportation Research Part B: Methodological}
  \bibinfo{volume}{40}, \bibinfo{pages}{335--350}.
\bibitem[{Fisher(1994)}]{fisher1994optimal}
\bibinfo{author}{Fisher, M.L.}, \bibinfo{year}{1994}.
\newblock \bibinfo{title}{Optimal solution of vehicle routing problems using
  minimum k-trees}.
\newblock \bibinfo{journal}{Operations research} \bibinfo{volume}{42},
  \bibinfo{pages}{626--642}.
\bibitem[{Franceschetti et~al.(2017)Franceschetti, Jabali and
  Laporte}]{franceschetti2017continuous}
\bibinfo{author}{Franceschetti, A.}, \bibinfo{author}{Jabali, O.},
  \bibinfo{author}{Laporte, G.}, \bibinfo{year}{2017}.
\newblock \bibinfo{title}{Continuous approximation models in freight
  distribution management}.
\newblock \bibinfo{journal}{Top} \bibinfo{volume}{25},
  \bibinfo{pages}{413--433}.
\bibitem[{Garn(2018)}]{Garn2018}
\bibinfo{author}{Garn, W.}, \bibinfo{year}{2018}.
\newblock \bibinfo{title}{Introduction to Management Science: Modelling,
  Optimisation and Probability}.
\newblock \bibinfo{publisher}{Smartana Ltd}.
\bibitem[{Garn(2020)}]{garn2020closed}
\bibinfo{author}{Garn, W.}, \bibinfo{year}{2020}.
\newblock \bibinfo{title}{Closed form distance formula for the balanced
  multiple travelling salesmen}.
\newblock \href{http://arxiv.org/abs/2001.07749}{{\tt arXiv:2001.07749}}.
\bibitem[{Gendreau et~al.(1999)Gendreau, Guertin, Potvin and
  Taillard}]{gendreau1999parallel}
\bibinfo{author}{Gendreau, M.}, \bibinfo{author}{Guertin, F.},
  \bibinfo{author}{Potvin, J.Y.}, \bibinfo{author}{Taillard, E.},
  \bibinfo{year}{1999}.
\newblock \bibinfo{title}{Parallel tabu search for real-time vehicle routing
  and dispatching}.
\newblock \bibinfo{journal}{Transportation science} \bibinfo{volume}{33},
  \bibinfo{pages}{381--390}.
\bibitem[{Gon{\c{c}}alves and Resende(2011)}]{Goncalves2011}
\bibinfo{author}{Gon{\c{c}}alves, J.F.}, \bibinfo{author}{Resende, M.G.},
  \bibinfo{year}{2011}.
\newblock \bibinfo{title}{{Biased random-key genetic algorithms for
  combinatorial optimization}}.
\newblock \bibinfo{journal}{Journal of Heuristics} \bibinfo{volume}{17},
  \bibinfo{pages}{487--525}.
\newblock \DOIprefix\doi{10.1007/s10732-010-9143-1}.
\bibitem[{Gouveia and Salazar-Gonz{\'a}lez(2010)}]{gouveia2010vehicle}
\bibinfo{author}{Gouveia, L.}, \bibinfo{author}{Salazar-Gonz{\'a}lez, J.J.},
  \bibinfo{year}{2010}.
\newblock \bibinfo{title}{On the vehicle routing problem with lower bound
  capacities}.
\newblock \bibinfo{journal}{Electronic Notes in Discrete Mathematics}
  \bibinfo{volume}{36}, \bibinfo{pages}{1001--1008}.
\bibitem[{Haghani and Jung(2005)}]{haghani2005dynamic}
\bibinfo{author}{Haghani, A.}, \bibinfo{author}{Jung, S.},
  \bibinfo{year}{2005}.
\newblock \bibinfo{title}{A dynamic vehicle routing problem with time-dependent
  travel times}.
\newblock \bibinfo{journal}{Computers \& operations research}
  \bibinfo{volume}{32}, \bibinfo{pages}{2959--2986}.
\bibitem[{James et~al.(2013)James, Witten, Hastie and
  Tibshirani}]{james2013introduction}
\bibinfo{author}{James, G.}, \bibinfo{author}{Witten, D.},
  \bibinfo{author}{Hastie, T.}, \bibinfo{author}{Tibshirani, R.},
  \bibinfo{year}{2013}.
\newblock \bibinfo{title}{An introduction to statistical learning}. volume
  \bibinfo{volume}{112}.
\newblock \bibinfo{publisher}{Springer}.
\bibitem[{Jaw et~al.(1986)Jaw, Odoni, Psaraftis and Wilson}]{jaw1986heuristic}
\bibinfo{author}{Jaw, J.J.}, \bibinfo{author}{Odoni, A.R.},
  \bibinfo{author}{Psaraftis, H.N.}, \bibinfo{author}{Wilson, N.H.},
  \bibinfo{year}{1986}.
\newblock \bibinfo{title}{A heuristic algorithm for the multi-vehicle advance
  request dial-a-ride problem with time windows}.
\newblock \bibinfo{journal}{Transportation Research Part B: Methodological}
  \bibinfo{volume}{20}, \bibinfo{pages}{243--257}.
\bibitem[{Kara and Bektas(2006)}]{kara2006integer}
\bibinfo{author}{Kara, I.}, \bibinfo{author}{Bektas, T.}, \bibinfo{year}{2006}.
\newblock \bibinfo{title}{Integer linear programming formulations of multiple
  salesman problems and its variations}.
\newblock \bibinfo{journal}{European Journal of Operational Research}
  \bibinfo{volume}{174}, \bibinfo{pages}{1449--1458}.
\bibitem[{Kirchler and Calvo(2013)}]{kirchler2013granular}
\bibinfo{author}{Kirchler, D.}, \bibinfo{author}{Calvo, R.W.},
  \bibinfo{year}{2013}.
\newblock \bibinfo{title}{A granular tabu search algorithm for the dial-a-ride
  problem}.
\newblock \bibinfo{journal}{Transportation Research Part B: Methodological}
  \bibinfo{volume}{56}, \bibinfo{pages}{120--135}.
\bibitem[{Kulak et~al.(2012)Kulak, Sahin and Taner}]{kulak2012joint}
\bibinfo{author}{Kulak, O.}, \bibinfo{author}{Sahin, Y.},
  \bibinfo{author}{Taner, M.E.}, \bibinfo{year}{2012}.
\newblock \bibinfo{title}{Joint order batching and picker routing in single and
  multiple-cross-aisle warehouses using cluster-based tabu search algorithms}.
\newblock \bibinfo{journal}{Flexible services and manufacturing journal}
  \bibinfo{volume}{24}, \bibinfo{pages}{52--80}.
\bibitem[{Lois and Ziliaskopoulos(2017)}]{lois2017online}
\bibinfo{author}{Lois, A.}, \bibinfo{author}{Ziliaskopoulos, A.},
  \bibinfo{year}{2017}.
\newblock \bibinfo{title}{Online algorithm for dynamic dial a ride problem and
  its metrics}.
\newblock \bibinfo{journal}{Transportation research procedia}
  \bibinfo{volume}{24}, \bibinfo{pages}{377--384}.
\bibitem[{Madsen et~al.(1995)Madsen, Ravn and Rygaard}]{madsen1995heuristic}
\bibinfo{author}{Madsen, O.B.}, \bibinfo{author}{Ravn, H.F.},
  \bibinfo{author}{Rygaard, J.M.}, \bibinfo{year}{1995}.
\newblock \bibinfo{title}{A heuristic algorithm for a dial-a-ride problem with
  time windows, multiple capacities, and multiple objectives}.
\newblock \bibinfo{journal}{Annals of operations Research}
  \bibinfo{volume}{60}, \bibinfo{pages}{193--208}.
\bibitem[{Martinez-Sykora and Bekta{\c{s}}(2015)}]{martinez2015transformations}
\bibinfo{author}{Martinez-Sykora, A.}, \bibinfo{author}{Bekta{\c{s}}, T.},
  \bibinfo{year}{2015}.
\newblock \bibinfo{title}{Transformations of node-balanced routing problems}.
\newblock \bibinfo{journal}{Naval Research Logistics (NRL)}
  \bibinfo{volume}{62}, \bibinfo{pages}{370--387}.
\bibitem[{Morais et~al.(2014)Morais, Mateus and Noronha}]{morais2014iterated}
\bibinfo{author}{Morais, V.W.}, \bibinfo{author}{Mateus, G.R.},
  \bibinfo{author}{Noronha, T.F.}, \bibinfo{year}{2014}.
\newblock \bibinfo{title}{Iterated local search heuristics for the vehicle
  routing problem with cross-docking}.
\newblock \bibinfo{journal}{Expert Systems with Applications}
  \bibinfo{volume}{41}, \bibinfo{pages}{7495--7506}.
\bibitem[{Navas(2017)}]{Navas2021}
\bibinfo{author}{Navas, M.}, \bibinfo{year}{2017}.
\newblock \bibinfo{title}{Taxi routes of mexico city, quito and more}.
\newblock \URLprefix
  \url{https://www.kaggle.com/mnavas/taxi-routes-for-mexico-city-and-quito/version/3}.
\bibitem[{Penna et~al.(2013)Penna, Subramanian and Ochi}]{Penna2013}
\bibinfo{author}{Penna, P.H.V.}, \bibinfo{author}{Subramanian, A.},
  \bibinfo{author}{Ochi, L.S.}, \bibinfo{year}{2013}.
\newblock \bibinfo{title}{{An iterated local search heuristic for the
  heterogeneous fleet vehicle routing problem}}.
\newblock \bibinfo{journal}{Journal of Heuristics} \bibinfo{volume}{19},
  \bibinfo{pages}{201--232}.
\newblock \URLprefix
  \url{https://link.springer.com/article/10.1007/s10732-011-9186-y},
  \DOIprefix\doi{10.1007/s10732-011-9186-y}.
\bibitem[{Pillac et~al.(2013)Pillac, Gendreau, Gu{\'e}ret and
  Medaglia}]{pillac2013review}
\bibinfo{author}{Pillac, V.}, \bibinfo{author}{Gendreau, M.},
  \bibinfo{author}{Gu{\'e}ret, C.}, \bibinfo{author}{Medaglia, A.L.},
  \bibinfo{year}{2013}.
\newblock \bibinfo{title}{A review of dynamic vehicle routing problems}.
\newblock \bibinfo{journal}{European Journal of Operational Research}
  \bibinfo{volume}{225}, \bibinfo{pages}{1--11}.
\bibitem[{Psaraftis(1980)}]{psaraftis1980dynamic}
\bibinfo{author}{Psaraftis, H.N.}, \bibinfo{year}{1980}.
\newblock \bibinfo{title}{A dynamic programming solution to the single vehicle
  many-to-many immediate request dial-a-ride problem}.
\newblock \bibinfo{journal}{Transportation Science} \bibinfo{volume}{14},
  \bibinfo{pages}{130--154}.
\bibitem[{Ratliff and Rosenthal(1983)}]{ratliff1983order}
\bibinfo{author}{Ratliff, H.D.}, \bibinfo{author}{Rosenthal, A.S.},
  \bibinfo{year}{1983}.
\newblock \bibinfo{title}{Order-picking in a rectangular warehouse: a solvable
  case of the traveling salesman problem}.
\newblock \bibinfo{journal}{Operations research} \bibinfo{volume}{31},
  \bibinfo{pages}{507--521}.
\bibitem[{Ruiz et~al.(2019)Ruiz, Soto-Mendoza, {Ruiz Barbosa} and
  Reyes}]{Ruiz2019}
\bibinfo{author}{Ruiz, E.}, \bibinfo{author}{Soto-Mendoza, V.},
  \bibinfo{author}{{Ruiz Barbosa}, A.E.}, \bibinfo{author}{Reyes, R.},
  \bibinfo{year}{2019}.
\newblock \bibinfo{title}{{Solving the open vehicle routing problem with
  capacity and distance constraints with a biased random key genetic
  algorithm}}.
\newblock \bibinfo{journal}{Computers and Industrial Engineering}
  \bibinfo{volume}{133}, \bibinfo{pages}{207--219}.
\newblock \DOIprefix\doi{10.1016/j.cie.2019.05.002}.
\bibitem[{Smolic-Rocak et~al.(2009)Smolic-Rocak, Bogdan, Kovacic and
  Petrovic}]{smolic2009time}
\bibinfo{author}{Smolic-Rocak, N.}, \bibinfo{author}{Bogdan, S.},
  \bibinfo{author}{Kovacic, Z.}, \bibinfo{author}{Petrovic, T.},
  \bibinfo{year}{2009}.
\newblock \bibinfo{title}{Time windows based dynamic routing in multi-agv
  systems}.
\newblock \bibinfo{journal}{IEEE Transactions on Automation Science and
  Engineering} \bibinfo{volume}{7}, \bibinfo{pages}{151--155}.
\bibitem[{Uchoa et~al.(2017)Uchoa, Pecin, Pessoa, Poggi, Vidal and
  Subramanian}]{uchoa2017new}
\bibinfo{author}{Uchoa, E.}, \bibinfo{author}{Pecin, D.},
  \bibinfo{author}{Pessoa, A.}, \bibinfo{author}{Poggi, M.},
  \bibinfo{author}{Vidal, T.}, \bibinfo{author}{Subramanian, A.},
  \bibinfo{year}{2017}.
\newblock \bibinfo{title}{New benchmark instances for the capacitated vehicle
  routing problem}.
\newblock \bibinfo{journal}{European Journal of Operational Research}
  \bibinfo{volume}{257}, \bibinfo{pages}{845--858}.
\bibitem[{Yang et~al.(2005)Yang, Hamedi and Haghani}]{yang2005online}
\bibinfo{author}{Yang, S.}, \bibinfo{author}{Hamedi, M.},
  \bibinfo{author}{Haghani, A.}, \bibinfo{year}{2005}.
\newblock \bibinfo{title}{Online dispatching and routing model for emergency
  vehicles with area coverage constraints}.
\newblock \bibinfo{journal}{Transportation Research Record}
  \bibinfo{volume}{1923}, \bibinfo{pages}{1--8}.

\end{thebibliography}
